\documentclass{amsart}

\usepackage{amssymb,amsfonts, amsmath, amsthm}
\usepackage[all,arc]{xy}
\usepackage{enumerate}
\usepackage{mathrsfs}
\usepackage{mathtools}
\usepackage[mathscr]{euscript}
\usepackage{array}
\usepackage{tikz}
\usepackage{tikz-cd}
\usepackage{textcomp}
\usetikzlibrary{decorations.pathmorphing,shapes}
\usepackage[pagebackref, colorlinks=true, linkcolor=blue, citecolor = red]{hyperref}
 \renewcommand*{\backref}[1]{}
 \renewcommand*{\backrefalt}[4]{({%
     \ifcase #1 Not cited.%
           \or On p.~#2%
           \else On pp.~#2%
     \fi%
     })}
\usepackage{todonotes}
\usepackage[nameinlink,capitalise,noabbrev]{cleveref}

\crefname{subsection}{Subsection}{Subsection}

\usetikzlibrary{patterns}

\usepackage{geometry}
\newcommand{\C}{\mathscr{C}}

\newcommand{\E}{\mathscr{E}}
\newcommand{\F}{\mathscr{F}}
\newcommand{\I}{\mathcal{I}}

\newcommand{\G}{\mathscr{G}}

\newcommand{\s}{\mathscr{S}}

\renewcommand{\O}{\mathscr{O}}
\renewcommand{\P}{\mathscr{P}}

\newcommand{\T}{\mathscr{T}}
\newcommand{\U}{\mathscr{U}}

\newcommand{\Q}{\mathscr{Q}}
\newcommand{\W}{\mathscr{W}}
\newcommand{\K}{\mathscr{K}}

\newcommand{\CSS}{\mathscr{C}\mathscr{S}\mathscr{S}}

\newcommand{\Map}{\mathrm{Map}}
\newcommand{\Hom}{\mathrm{Hom}}
\newcommand{\comma}{,}
\newcommand{\leb}{[}
\newcommand{\reb}{]}
\newcommand{\set}{\mathscr{S}\text{et}}
\newcommand{\cat}{\mathscr{C}\text{at}}

\newcommand{\Kan}{\mathscr{K}\text{an}}
\newcommand{\sSet}{\text{s}\mathscr{S}\text{et}}

\newcommand{\ssSet}{\text{ss}\mathscr{S}\text{et}}
\newcommand{\Sub}{\mathscr{S}\text{ub}}

\newcommand{\colim}{\mathrm{colim}}

\newcommand{\Ext}{\mathscr{E}\text{xt}}
\newcommand{\Comp}{\mathscr{C}\text{omp}}

\newcommand{\Los}{{\L}o{\'s}}

\newcommand{\citerefive}[1]{\cite[#1]{rezk2010toposes}}

\newcommand{\citenno}[1]{\cite[#1]{rasekh2021nno}}

\newcommand{\comsq}[8]{
  \begin{tikzcd}[row sep=0.5in, column sep=0.5in]
    #1 \arrow[r, "#5"] \arrow[d, "#6"']
    \pgfmatrixnextcell #2 \arrow[d, "#7"] \\
    #3 \arrow[r, "#8"]
    \pgfmatrixnextcell #4
  \end{tikzcd}
}

\newcommand{\pbsq}[8]{
  \begin{tikzcd}[row sep=0.5in, column sep=0.5in]
    #1 \arrow[r, "#5"] \arrow[d, "#6"'] \arrow[dr, phantom, "\ulcorner", very near start]
    \pgfmatrixnextcell #2 \arrow[d, "#7"] \\
    #3 \arrow[r, "#8"']
    \pgfmatrixnextcell #4
  \end{tikzcd}
}

\newtheorem{theone}{Theorem}[section]
\newtheorem*{thetwo}{Theorem}
\newtheorem{lemone}[theone]{Lemma}
\newtheorem{propone}[theone]{Proposition}
\newtheorem{corone}[theone]{Corollary}

\theoremstyle{definition}
\newtheorem{defone}[theone]{Definition}
\newtheorem{exone}[theone]{Example}
\newtheorem*{extwo}{Example}

\newtheorem{notone}[theone]{Notation}

\newtheorem{convone}[theone]{Convention}

\theoremstyle{remark}
\newtheorem{remone}[theone]{Remark}

\makeatletter
\def\@seccntformat#1{%
  \expandafter\ifx\csname c@#1\endcsname\c@section\else
  \csname the#1\endcsname\quad
  \fi}
\makeatother

\title{Filter Quotients and Non-Presentable $(\infty,1)$-Toposes}

\author{Nima Rasekh}

\date{April 2021}

\address{{\'E}cole Polytechnique F{\'e}d{\'e}rale de Lausanne, SV BMI UPHESS, Station 8, CH-1015 Lausanne, Switzerland}

\email{nima.rasekh@epfl.ch}

\begin{document}

\begin{abstract}
We define filter quotients of $(\infty,1)$-categories and prove that filter quotients preserve the structure of an
elementary $(\infty,1)$-topos and in particular lift the filter quotient
of the underlying elementary topos.
We then specialize to the case of filter products of $(\infty,1)$-categories and 
prove a characterization theorem for equivalences in a filter product.
\par 
Then we use filter products to construct a large class of elementary $(\infty,1)$-toposes that 
are not Grothendieck $(\infty,1)$-toposes. Moreover, we give one detailed example for the interested reader who would like to 
see how we can construct such an $(\infty,1)$-category, but would prefer to avoid the technicalities regarding filters.
\end{abstract}

\maketitle
\addtocontents{toc}{\protect\setcounter{tocdepth}{1}}

\tableofcontents

 \section{Introduction} \label{Sec:Introduction}
 
 \subsection{The Set-Theoretical Foundations of Algebraic Topology}
 The study of {\it algebraic topology} has historically strongly relied on a set-theoretical foundation.
 Many classical results that we would now associate to the field of algebraic topology were proven using some notion of
 {\it topological spaces} (manifolds, complexes, ...), which are quite literally sets along with extra data
 \cite{poincare1895analysissitus, hopf1931fibration, whitehead1949combinatoriali, whitehead1949combinatorialii}.
 
 Later, an alternative approach towards algebraic topological computations emerged: {\it simplicial sets} \cite{eilenbergzilber1950simplicial, kan1957complexes, kan1958homotopygroups}. However, early results about simplicial sets relied again on having set-theoretical foundations. It was later observed that both topological spaces and simplicial sets carry {\it model structures} and they are in fact equivalent \cite{quillen1967modelcats}, which formalized the intuition that both can be used as a foundation for algebraic topology. 
 
 The trend established by these structures has continued to dominate the study of algebraic topology until today. 
 For example, while computing homotopy groups, and in particular homotopy groups of spheres, we can observe certain stability phenomena, that 
 lead to {\it stable homotopy theory} and the study of {\it spectra} \cite{lima1959spectrum, adams1995stable}. There are various ways of characterizing spectra, but they all involve a collection of spaces with certain relations \cite{ekmm1995stable, hss2000symmetricspectra, mandellmay2002orthogonal}.
 Hence, the study of stable homotopy theory continued the dependence on set-theoretical foundations.
 
 \subsection{Elementary Toposes as a Foundation}
 In that same time frame a very rich theory of mathematical foundations apart from the standard set theory familiar to algebraic topologists emerged. There was an extensive analysis of {\it first-order (elementary) languages} (mathematical settings without any set-theoretical assumptions) in the context of {\it model theory} \cite{godel1930completeness, tarski1956semantics, changkeisler1990modeltheory}. Parallel to that development, there also emerged an alternative to set theory and its model theoretical approach via the notion of {\it type theory} \cite{russell1908typetheory, church1940typetheory}.
 At that time, such foundational studies seemed very distant from algebraic topological concerns.
 
 A first step towards observing meaningful connections between those two was the development of {\it category theory} \cite{eilenbergmaclane1945categories}.
 Although it first emerged in the context of algebraic topology with the goal of studying naturality \cite{eilenbergsteenrod1945axioms}, it very quickly became foundational in other branches of mathematics, and in particular algebraic geometry, where the school of Grothendieck and Bourbaki realized the centrality of the {\it Grothendieck topos} in order to further study sheaf theory \cite{sga1972tome1, sga1972tome2, sga1972tome3}.
 It was then realized by Lawvere and Tierney, building on previous work \cite{lawvere1964elementarysets}, that this definition can be generalized to an {\it elementary topos}, which is a category that can itself serve as a foundation to mathematics and be another alternative to set theory \cite{tierney1973elementarytopos}.
 
 Using an elementary topos (some times with additional assumptions) they, and other authors, construct many interesting foundations for mathematics:
 \begin{enumerate}
 	\item {\bf Models of set theory:} It was proven that {\it well-pointed Boolean elementary toposes with natural number object} correspond to categories with objects sets that satisfy the {\it Zermelo-Fraenkel axioms} \cite{maclanemoerdijk1994topos}. This implies that we can construct models of set-theory by constructing certain elementary toposes. As a direct application we can prove that the {\it continuum hypothesis} is independent of the other axioms of ZF by constructing appropriate categories \cite{tierney1972elementarycontinuum}. Hence, giving an alternative proof to the original one by Cohen, that used the {\it forcing technique} \cite{cohen1963forcing}. 
 	
 	\item {\bf Filter Quotient:} One benefit of using a categorical approach to foundations is the fact that there are many ways to construct new categories and hence many new foundations. For examples, motivated by the study of ultra products in model theory \cite{changkeisler1990modeltheory}, there is a filter quotient construction for elementary toposes \cite{volger1975filter}. Moreover, Adelman and Johnstone study which properties of the elementary topos will transfer along the filter quotient \cite{adelmanjohnstone1982serreclasses}. One implication is that we can apply the filter construction to construct very unusual models of set theory, which even lack countable coproducts \cite[Example A2.1.13]{johnstone2002elephanti}. 
 	
 	\item {\bf Non-Standard Analysis} Applying the ultra product construction to certain sheaves on real numbers, we can develop a non-standard analysis internal to the topos giving us new ways to tackle the concept of infinity and convergence \cite{palmgren1997nonstandardanalysis}.
 	 \item {\bf Intuitionistic Logic:} As an example related to the previous one, in 1927 Brouwer proved that all real valued functions are continuous \cite{brouwer1927allcontinuous}! The way to make sense of this fact is that Brouwer used {\it intuitionistic logic}, which is constructive in nature. Using topos theory, we can in fact construct a topos where this holds \cite{maclanemoerdijk1994topos}.
 \end{enumerate}
 The key fact here is that the definition of an elementary topos is in fact {\it elementary}, or {\it first-order},  meaning the definition does not rely on any set-theoretical assumptions. In particular, it manages to avoid any infinite constructions that would necessitate first establishing necessary set-theoretical axioms that allows us to handle the concept of infinity, such as infinite limits and colimits.
 
 We can use the elementary nature of the internal language of an elementary topos to prove statements inside the topos using appropriate syntactic rules, which has been made precise by studying its connection with type theory.
 We can prove that every elementary topos is a model for certain {\it higher-order type theories} and so we can prove statements inside the elementary topos by constructing the appropriate type inside the type theory \cite{lambekscott1988higherorderlogic}. 
 This equivalence allows us to effectively treat a topos the way we would treat the category of sets. For interesting examples illustrating how this is implemented in practice see \cite[D5.1]{johnstone2002elephantsii}.
  
 \subsection{Why have different foundations for algebraic topology?}
 At a first glance, we might disregard such discussion of foundations as not pertinent to significant questions in algebraic topology. 
 However, a second glance does suggest some interesting examples:
 \begin{enumerate}
 	\item {\bf Natural number object:} In algebraic topology we instinctively order the homotopy groups by the set $\mathbb{N}= \{0,1,2,...\}$. 
 	By extension, spectra are then ordered by the integers, which is simply the group completion. 
 	This comes from the fact that $\mathbb{N}$ is the {\it natural number object} inside the standard category of sets, which, as discussed before, has been the foundation of algebraic topology until now. 
 	
 	There are examples of elementary toposes where the natural number object is not standard \cite[D5.1.7]{johnstone2002elephantsii}, which opens the possibility of defining non-standard homotopy groups indexed by the non-standard natural number object. By extension this would also change how we index a spectrum and so also influences the stable homotopy theory in this setting. 
 	
 	However, in order to seriously tackle such questions, we would first need to properly develop algebraic topology in such a general setting.
 	\item {\bf Homotopy type theory:} In the $70$'s Martin-L{\"o}f introduced what is now called {\it Martin-L{\"o}f type theory}, introducing the notion of a {\it identity type} \cite{martinlof1975inttypetheories, martinlof1975inttypetheoriespredicative}, giving us a starting point for introducing homotopical thinking into type theory. 
 	This was brought to fruition by Voevodsky who introduced {\it homotopy type theory} or {\it univalent foundations} and used that to prove classical algebraic topological results in this setting \cite{hottbook2013}. In the long run Voevodsky was hoping to formalize advanced homotopy theory in a proof assistant in order to prevent the preponderance of mistakes in a theory that is constantly getting more complicated \cite{voevodsky2014origins}. One side benefit of univalent foundations is that it proves results from algebraic topology in a manner independent from set theory. However, this also poses a challenge. Concretely whenever a proof in algebraic topology explicitly relies on set theoretical assumptions then it cannot be translated immediately into a proof in the homotopy type theoretical setting and thus cannot be formalized. 
 	
 	Hence, in order to further the formalization of homotopy theory we need to be able to develop algebraic topological constructions in different foundations. 
 	
    \item {\bf Non-standard models of spaces:} Even if we were solely focused on algebraic topology with standard set-theoretical foundations then taking a broader view can have immense value. For example, while trying to study the homotopy theory of schemes, Voevodsky developed {\it motivic homotopy theory} \cite{voevodsky1998motivicorigin}. It has found many applications in algebraic geometry, such as solving the {\it Bloch-Kato conjecture} 
    \cite{suslinvoevodsky2000blochkato}, however, it has also helped us further compute classical homotopy groups of spheres, as can be witnessed in the work of Isaksen \cite{isaksen2019stablecomputation}. Somewhat ironically, the additional complexity of motivic spheres enables certain computations that would have otherwise been far more challenging.
    This observation naturally leads to the question whether there are other interesting alternatives to spaces which can help us further our understanding of the homotopy groups of spheres, which again necessitates studying algebraic topology in different foundations.
    
 	\item {\bf Constructing Kan model structure:} The {\it Kan model structure} on simplicial sets has been one of the most important model structures in homotopy theory. Until recently it was assumed that any construction of the Kan model structure would require all axioms of set theory \cite{quillen1967modelcats, goerssjardine1999simplicialhomotopytheory, hirschhorn2003modelcategories}. 
 	
 	By taking a constructive approach to mathematics, Gambino-Sattler-Szumilo \cite{gss2019constructive} and Henry \cite{henry2019constructive} construct a model structure on simplicial sets without assuming either the law of excluded middle or the axiom of choice, giving us a possibility to do algebraic topology in a much broader setting. 
 	
 	In order to determine how much further this approach can be taken and which foundations would permit us to construct a version of the Kan model structure on simplicial objects we need to first better understand algebraic topology in alternative foundations.
 	
 \end{enumerate}
 
 This leaves us with many important questions:
 
 \begin{enumerate}
 	\item We can understand sets collectively by studying its category. What is a proper framework to study topological spaces?
 	\item Can we find analogous definitions to Grothendieck toposes and elementary toposes?
 	\item Can we prove results analogous to what we have observed for elementary toposes in such a setting?
 \end{enumerate}

\subsection{$\infty$-Toposes}
Fortunately many of these questions have already been answered. Concretely, we now consider {\it $(\infty,1)$-categories} \cite{bergner2010survey} (and sometimes model categories \cite{quillen1967modelcats, hovey1999modelcategories}) as the appropriate generalizations of categories, suitable for the study of homotopy theories, and in particular the homotopy theory of topological spaces.
We also have a working generalization of Grothendieck toposes to this $(\infty,1)$-categorical setting, known as {\it $\infty$-toposes} \cite{lurie2009htt} (also called {\it model topos} \cite{rezk2010toposes} in the model categorical setting). Henceforth, we shall call them {\it Grothendieck $(\infty,1)$-toposes} to avoid confusion with other terminology. 

Doing algebraic topology inside a Grothendieck $(\infty,1)$-topos is already a great step in the right direction. In particular, it is already established that it models homotopy type theory \cite{shulman2019inftytoposunivalent} and we can use that setting to prove classical results in a much more set theory independent manner, as can be seen for the examples of the {\it Blakers-Massey theorem} \cite{abfj2017blakersmassey}, {\it Goodwillie calculus} \cite{abfj2018goodwillie} and {\it localization theory} \cite{vergura2019localization}. 
Moreover, Grothendieck higher topos theory has been used as a foundation to develop a lot of {\it derived} (and {\it spectral}) {\it geometry} \cite{lurie2009htt, lurie2017ha, lurie2018sag}.

Following the previous trends, the next step should then be to generalize Grothendieck $(\infty,1)$-toposes to {\it elementary $(\infty,1)$-toposes} and to prove the desired results. A definition does already exist \cite{rasekh2018elementarytopos}, but most desired results remain unsolved.
Concretely, here are certain results and remaining work: 
\begin{enumerate}
	\item We do know that every Grothendieck $(\infty,1)$-topos is in fact an elementary $(\infty,1)$-topos \cite{rasekh2018elementarytopos}. 
	Moreover, we do have some examples of elementary $(\infty,1)$-toposes that are not Grothendieck \cite{monaco2021eht}, however, those examples are sub-categories of the category of spaces and so do not differ from a foundational perspective (for example have the same natural number object).
	\item We can recover some classical algebraic topological constructions in this setting (such as truncations, Blakers-Massey theorem) \cite{rasekh2018truncations}, 
	however, certain results have remained unsettled (for example homotopy groups of spheres, spectral sequence computations, stabilizations).
	\item We do not generally understand their relation with homotopy type theory, except in the cases already covered by Grothendieck $(\infty,1)$-toposes \cite{shulman2019inftytoposunivalent}.
\end{enumerate}

\subsection{Filter Quotients of $(\infty,1)$-Categories}
 The goal of this paper is to take a further step towards better understanding higher topos theory by studying {\it filter quotients}. 
  Concretely, we will define a filter quotient for $(\infty,1)$-categories and then prove that elementary $(\infty,1)$-toposes are closed under this construction (see \cref{Subsec:Main Results} for more detailed description of the main results). 
  This will then in particular lead to a large class of examples of elementary $(\infty,1)$-toposes where the underlying logic is quite non-standard (we will carefully analyze one example in \cref{Subsec:An Example}).
  
  We are here focusing on the benefits of the filter construction to higher topos theory, as our primary goal is to construct non-trivial examples. However, we do expect interesting applications in the opposite direction, and concretely {\it chromatic homotopy theory}.
  In chromatic homotopy theory we usually break down spectra in two consecutive steps. First, we choose a prime $p$ and then 
  we choose a height $n$ (corresponding to the various {\it Morava $K$-theories}). As the heights form a natural filtration
  we can then study the limiting behavior as the height increases, as long as we fix a prime.
  \par 
  On the other hand, the various prime numbers do not give us a filtration and thus there is no immediate way to study the limiting behavior 
  as we increase the prime numbers. The insight of Barthel, Schlank and Stapleton \cite{bss2020chromatic} was to instead study the ultra product 
  of the various $p$-local subcategories and show that the ultra product is equivalent to a category that comes from algebraic geometry, at the same 
  time giving meaning and proving the statement ``chromatic homotopy theory is asymptotically algebraic at a fixed height".
  \par 
  However, when setting up the theory of ultra products, they observe that general ultra product $(\infty,1)$-categories are not well-behaved 
  \cite[Example 3.19]{bss2020chromatic}, in the sense that the ultra product construction does not preserve presentability. This presents a serious challenge, as most of the results in stable homotopy theory have been proven for presentable stable $(\infty,1)$-categories. 
  Hence, they are forced to use an alternative construction,  the {\it compactly generated ultra product}, that does in fact leads to a presentable $(\infty,1)$-category, however, does not precisely coincide with the original construction anymore \cite[Section 3.3]{bss2020chromatic}. 
  \par 
  The power of elementary $(\infty,1)$-topos theory is that its definition does not rely on presentability. Hence any result proven about them simply transfers to the ultra product and we do not need to make any adjustments. This opens the possibility of studying the ultra product categories that arise in the work of \cite{bss2020chromatic}, and are of interest to chromatic homotopy theorists, directly. 
  
 \subsection{Main Results} \label{Subsec:Main Results}
 The main result can be summarized as the following theorems:
 
 First, we introduce a general method for constructing a new $(\infty,1)$-category, namely the {\it filter quotient}:
 
 \begin{thetwo}
  Let $\C$ be a finitely complete $(\infty,1)$-category, which can be:
  \begin{enumerate}
   \item a complete Segal space,
   \item a quasi-category,
   \item a Kan enriched category
  \end{enumerate}
  and $\Phi$ a filter (\cref{Def:Filter}) of subobjects (\cref{Ex:Filter Subob}).
  \par 
  Then there exists a finitely complete $(\infty,1)$-category $\C_\Phi$ (\cref{Prop:E Phi simp Cat}, 
  \cref{Def:Filter Quotient CSS}, \cref{The:Filter Quotient QCat})
  along with a functor $P_\Phi: \C \to \C_\Phi$ (\cref{Def:P Phi CSS}, \cref{Def:P Phi Kan Enrich}).
 \end{thetwo}

 Next we show the construction is well behaved from the perspective of elementary topos theory:
 
 \begin{thetwo}
  Let $\E$ be an elementary $(\infty,1)$-topos and $\Phi$ a filter of subobjects. 
  \begin{enumerate}
   \item (\cref{The E Phi EHT}) 
   The functor $P_\Phi$ preserves all defining conditions of an elementary $(\infty,1)$-topos:
   \begin{enumerate}
    \item finite (co)limits
    \item local Cartesian closure
    \item subobject classifiers
    \item universes
   \end{enumerate}
   and so $\E_\Phi$ is an elementary $(\infty,1)$-topos as well.
   \item (\cref{The:E Phi Tau egal E Tau Phi}) We have an equivalence of underlying elementary toposes $\tau_0(\E_\Phi) \simeq \tau_0(\E)_\Phi$.
   \item (\cref{The:P Phi preserves NNO}) So, $P_\Phi$ also preserves natural number objects in $\E_\Phi$.
   \item (\cref{Cor EHT not GHT}) $\E_\Phi$ is not a Grothendieck $(\infty,1)$-topos 
   (\cref{Def:Groth Infty Topos}) if $\tau_0(\E)_\Phi$ is not a Grothendieck topos. 
  \end{enumerate}
 \end{thetwo}
 
 This gives us a very effective recipe to construct non-presentable Grothendieck $(\infty,1)$-toposes, 
 namely constructing filter quotients $\E_\Phi$ such that the underlying elementary topos $\tau_0(\E)_\Phi$ is not presentable.
 \par 
 In order to find such examples we next focus on a specific class of filter quotients, namely the {\it filter product}.
 We first show that in some ways they behave similar to the classical analogue:
 
 \begin{thetwo}
  Let $\C$ be a finitely complete $(\infty,1)$-category such that the final object has only two subobjects, $I$ be a set and $\Phi$ 
  a filter on $P(I)$. Then we can construct the filter product $\prod_\Phi \C$ and we have the following \Los \ type results:
  \begin{enumerate}
   \item (\cref{The:Los Two maps Equiv}, \cref{Cor:Los Two maps Equiv}) 
   Two maps $(f_i)_{i \in I}, (g_i)_{i \in I}$ are equivalent if and only if 
   $$\{ i \in I : f_i \simeq g_i \} \in \Phi $$ 
   \item (\cref{The:Los Maps Equiv}, \cref{Cor:Los Maps Equiv}) A map $(f_i)_{i \in I}$ is an equivalence if and only if 
   $$\{ i \in I : f_i \text{ is an equivalence} \} \in \Phi $$ 
  \end{enumerate}
 \end{thetwo}
 
 Finally, we use the filter product to give a large of examples of elementary $(\infty,1)$-toposes that are not Grothendieck $(\infty,1)$-toposes:
 
 \begin{extwo}
  (\cref{Ex:KanS Phi not GHT})
  There is a class of sets $I$ and filter $\Phi$ such that the filter product $\prod_\Phi \set$ is not a Grothendieck topos 
  and thus $\prod_\Phi \Kan$ is an elementary $(\infty,1)$-topos that is not a Grothendieck $(\infty,1)$-topos.
 \end{extwo}

 \subsection{Outline}
 \cref{Sec:Review of Concepts} gives a review of some of the concepts we need later on (see also \cref{Subsec:Background}
 for an overview of the necessary background and \cref{Subsec:Notation} for notation).
 In particular, we review several important features of complete Segal spaces (\cref{Subsec:Theory of Infty Categories})
 that play a crucial role in \cref{Sec:The Filter Quotient}.
 \par
 \cref{Sec:The Filter Quotient} breaks down into two subsections.
 In the first subsection (\cref{Subsec:The Filter Quotient}) we define the filter quotient. We give two separate definitions for 
 Kan enriched categories and for complete Segal spaces and then prove these are equivalent. 
 In the second subsection (\cref{Subsec:Filter Quotient EHT}) we prove that the filter quotient construction preserves 
 all the properties of an elementary $(\infty,1)$-topos.
 \par
 \cref{Sec:Filter Products} focuses on a special filter quotient: {\it filter products}. In 
 \cref{Subsec:Los Theorem for Equivalences} we study general features of filter products and prove some 
 results along the lines of \Los's theorem. In \cref{Subsec:An Example} we finally use filter products to 
 give examples of elementary $(\infty,1)$-toposes that are not Grothendieck $(\infty,1)$-toposes. 
 We will present one example in great detail, for the benefit of any reader who is only interested in understanding one example and wants
 to skip the technicalities about filters quotients.
 \par
 Finally, \cref{Sec:Future} mentions some possible future directions. 
 
 \subsection{Background} \label{Subsec:Background}
 We will assume the reader has already some familiarity with the following concepts:
 \begin{enumerate}
  \item The theory of $(\infty,1)$-categories, in particular the different models and how they are related.
  \item The theory of elementary $(\infty,1)$-toposes.
 \end{enumerate}
 We only give a quick overview in \cref{Sec:Review of Concepts} and refer the reader to the appropriate sources.
 \par 
 On the other hand, the following concepts require no familiarity and everything we need is covered in 
 \cref{Sec:Review of Concepts}:
 \begin{enumerate}
  \item Grothendieck $(\infty,1)$-topos theory,
  \item Filters.
 \end{enumerate}
 
  \subsection{Notation} \label{Subsec:Notation}
  We are using various models of $(\infty,1)$-categories with various levels of strictness. In order to avoid any confusion we 
  will use the following conventions:
  \begin{enumerate}
   \item $\set$, $\sSet$ and $\ssSet$ are the $1$-categories of sets, simplicial sets and bisimplicial sets (also called simplicial spaces), respectively.
   \item The notation $\colim$ refers to a colimit in a strict $1$-category.
   \item The notation $\cong$ refers to an isomorphism of objects in a $1$-category.
   \item We will use Kan complexes as our preferred model for the ``homotopy theory of spaces" and therefore will 
   avoid using the vague term space throughout.
   \item As we use various models of $(\infty,1)$-categories we use a (rather arbitrary) convention to help the reader 
   distinguish the various notations.
   \begin{itemize}
    \item[(I)] $\C$: Arbitrary $(\infty,1)$-Category
    \item[(II)] $\E$: Arbitrary elementary $(\infty,1)$-topos
    \item[(III)] $\G$: Arbitrary Grothendieck $(\infty,1)$-topos
    \item[(IV)] $\Q$: Quasi-Category
    \item[(V)] $\W$: Complete Segal space
    \item[(VI)] $\K$: Kan enriched category
   \end{itemize}

  \end{enumerate}
 
  \subsection{Acknowledgments}
  I want to thank Asaf Horev for pointing to the connection with the work by Barthel, Schlank and Stapleton \cite{bss2020chromatic}.
  I also want to thank Peter Lumsdaine for making me aware of a more general form of \Los's theorem that applies to non-ultrafilter,
  which resulted in the material in \cref{Subsec:Los Theorem for Equivalences}.
  Finally, I want to thank the referee for many helpful suggestions which in particular has resulted in a much better introduction and clarified several important proofs.
 
 \section{Review of Concepts} \label{Sec:Review of Concepts}
 We will give a minimal overview over the specific concepts required later on and refer the reader to the appropriate sources 
 for more details. 
 
 \subsection{The Theory of \texorpdfstring{$(\infty,1)$}{(oo,1)}-Categories} \label{Subsec:Theory of Infty Categories}
 An $(\infty,1)$-category is a general idea of a weak category or homotopical category i.e. 
 a category with a weak notion of composition and associativity. 
 Unlike classical category theory, making such an intuition into a precise definition has been quite challenging. Hence, the intuition of $(\infty,1)$-categories has been captured by various {\it models}, some of the most important ones are 
 {\it quasi-categories} \cite{boardmanvogt1973qcats},
 {\it Kan enriched categories} \cite{dwyerkan1980calculatingsimplocalizations, dwyerkan1980simplocalization}, 
 and {\it complete Segal spaces} \cite{rezk2001css}.

 Each of these models has their own model structures, namely {\it Joyal} \cite{joyal2008notes}, {\it Bergner} \cite{bergner2007bergnermodelcat} and {\it Rezk} \cite{rezk2001css} model structure respectively.
 Moreover, they are all Quillen equivalent via various Quillen equivalences. Concretely we have the following diagram of Quillen equivalences 
 
 \begin{center}
 	\begin{tikzcd}[row sep=0.5in, column sep=0.5in]
 		\sSet^{Joyal} \arrow[r, shift left=1.4, "p_1^*", "\sim"']  & 
 		\ssSet^{Rezk} \arrow[r, shift left=1.4, "t_!", "\sim"'] \arrow[l, shift left=1.4, "i_1^*"] & 
 		\sSet^{Joyal} \arrow[r, shift left=1.4, "\mathfrak{C} \leb - \reb", "\sim"'] \arrow[l, shift left=1.4, "t^!"] & 
 		(\cat_\Delta)^{Bergner} \arrow[l, shift left=1.4, "N_\Delta"]
 	\end{tikzcd}.
 \end{center}
 
 For the Quillen equivalence ($\mathfrak{C}[-]$, $N_\Delta$) see \cite{lurie2009htt} 
 and for the other two Quillen equivalences see \cite{joyaltierney2007qcatvssegal}.
 Moreover, for an overview about these and other models of $(\infty,1)$-categories see \cite{bergner2010survey, bergner2018book}.
 
 There are certain aspects of these model structures and their Quillen equivalences that will play an important role in the coming sections and hence will be quickly reviewed here. 
 First of all we benefit from the fact that right Quillen functors preserve fibrant objects. 
  
  \begin{remone} \label{Rem:NDelta sCat is CSS}
 	The functors $N_\Delta$ and $t^!$ are Quillen right adjoints and thus for any Kan enriched category $\K$, 
 	$N_\Delta(\K)$ is a quasi-category and $t^!N_\Delta(\K)$ is a complete Segal space. 
 \end{remone}
 
 Next, all three models have a notion of {\it mapping Kan complex}. For a Kan enriched category it is given as part of a data. 
 We will also need the mapping space of a complete Segal space.
 
 \begin{defone} \label{Def:CSS Mapping Space}
 	Let $\W$ be a complete Segal space.
 	For two objects $X,Y$ (i.e. points in $\W_0$) we define the {\it mapping Kan complex} as the pullback
 	\begin{center}
 		\pbsq{\Map_\W(X,Y)}{\W_1}{*}{\W_0 \times \W_0}{}{}{(s \comma t)}{ ( X \comma Y)}.
 	\end{center}
 \end{defone}
 
 Moreover, the Quillen equivalences relate the mapping spaces in an appropriate manner.
  
  \begin{remone} \label{Rem:Mapping Spaces sCat to CSS}
 	If we have a Kan enriched category $\K$ then the complete Segal space $t^!N_\Delta(\K)$ has the 
 	same objects as $\K$ and we have an equivalence of mapping Kan complexes
 	$$\Map_\K(X,Y) \simeq \Map_{t^!N_\Delta\K}(X,Y).$$
 \end{remone}
 
 Next, we observe that all three models have a notion of {\it underlying $(\infty,1)$-groupoid}, or {\it core}.
 Following the convention of classical category theory, for Kan enriched categories it is the subcategory with the same objects and morphisms equivalences. For the other two models we have the following definitions.
 
 \begin{defone} \label{Def:Core of QCat}
 	Let $\Q$ be a quasi-category. Then we call the maximal Kan complex inside $\Q$ the {\it core of $\Q$}, 
 	denoted $\Q^{core}$.
 \end{defone}
  
  \begin{defone} \label{Def:Core of CSS}
 	Let $\W$ be a complete Segal space. Then we define the underlying $(\infty,1)$-groupoid $\W$ as 
 	$\W^{core} = \W_0$.
 	Hence,  we also call $\W_0$ the {\it Kan complex of objects}. 
 \end{defone}
 
 \begin{remone}
 	Given that a complete Segal space $\W$ is itself a simplicial space, we might have expected the underlying $(\infty,1_)$-groupoid to be a simplicial space as well. 
 	However, by the completeness condition, a complete Segal space $\W$ is an $(\infty,1)$-groupoid if and only if it is {\it homotopically constant}, meaning that for any map $\delta:[n] \to [m]$ the corresponding map of Kan complexes $\delta^*:\W_n \to \W_m$ is a Kan equivalence. 
 	
 	In other words, $(\infty,1)$-groupoids are precisely the Reedy fibrant simplicial spaces that are weakly equivalent to constant simplicial spaces. 
 	Hence, up to equivalence such complete Segal spaces are characterized by their $0$-space, and there is no loss of information. 
 	
    On the other hand, defining the core of a complete Segal space as a simplicial set has the additional benefit that for a quasi-category $\Q$
 	$$(t^!\Q)^{core} = \Q^{core},$$
 	meaning the definition of the core for complete Segal spaces coincides with the definition for quasi-categories.
 \end{remone}
 
 Using the core we can construct fibrant objects out of left Quillen functors.
 
   \begin{remone} \label{Rem:Gamma of QCats}
 	For a given quasi-category $\Q$, $p_1^*(\Q)$ is not a complete Segal space (unless $\Q$ has 
 	no non-trivial automorphisms). So we need an alternative construction. Thus we introduce $\Gamma(\Q)$, defined 
 	as the simplicial space
 	$$\Gamma(\Q)_n = (\Q^{\Delta^n})^{core}$$ 
 	and observe that it is a complete Segal space with the property that $\Gamma(\Q)_{n0} = \Q_n$
 	(for more details see \cite{joyaltierney2007qcatvssegal}).
 \end{remone}
 
  There is one $(\infty,1)$-category that will play an important role later on.
  
 \begin{notone}
 	We denote by $\Kan$ the Kan-enriched category of Kan complexes. 
 	On the other hand, we use the notation $\s$ for the quasi-category and the complete Segal space of Kan complexes.
 	Note, in particular, $\s = N_\Delta(\Kan)$.
 \end{notone}

 Up until now we have not taken a model-independent approach to our treatment of $(\infty,1)$-categories. 
 There are now in fact various ways to develop $(\infty,1)$-categories model-independently, either via $2$-categories \cite{riehlverity2017inftycosmos} or via homotopy type theory \cite{riehlshulman2017rezktypes}. However, we will not work at this level of generality and rather set the following convention.
 
 \begin{convone}
 	The term {\it $(\infty,1)$-category} henceforth refers to quasi-categories, Kan enriched categories or complete Segal spaces. 
 \end{convone}
 
 This means, whenever a definition is well-defined for all three models, we will apply it to $(\infty,1)$-categories. 
 This convention will in particular apply to objects, morphism and the core.
 Here is another important example of a model independent definition.
 
 \begin{defone} \label{Def:LCCC}
 	Let $\C$ be an $(\infty,1)$-category with finite limits. Then $\C$ is {\it locally Cartesian closed} if for every morphism $f: x \to y$ the pullback functor 
 	$$f^*: \C_{/y} \to \C_{/x}$$
 	has a right adjoint $f_*: \C_{/x} \to \C_{/y}$.
 \end{defone}
 
 We will see further model independent definitions in the next sections.
 
  We also need some understanding of the basics of complete Segal spaces (as given in \cite{rezk2001css}), as well as over-categories and limits in the context of complete Segal spaces \cite{rasekh2017left}. 
  
   \begin{notone} \label{Not:Free Arrow CSS}
  	In the context of complete Segal spaces we denote the free arrow as $F(1)$, following notation in \cite{rezk2001css}.
  	Similarly, $F(0)$ denotes the complete Segal space with one object. In particular, both of these are simplicial spaces.
  \end{notone}

  \begin{defone} \label{Def:Over Cat CSS}
   Let $\W$ be a complete Segal space and $x$ an object. Then we can define the slice CSS $\W_{/x}$ as the pullback 
   \begin{center}
    \pbsq{\W_{/x}}{\W^{F(1)}}{F(0)}{\W}{}{}{t}{\{x\}}
   \end{center}
  \end{defone}
  
  \begin{defone}\label{Def:CSS of Cones}
   More generally, let $I$ be a simplicial space and $G: I \to \W$ be a fixed diagram in the CSS $\W$. Then 
   we define the {\it complete Segal space of cones} as
   $$\W_{/G} = F(0) \underset{\W^I}{\times} \W^{I \times F(1)} \underset{\W^I}{\times} \W$$
  \end{defone}
  
  Moreover, we have the following results about the cones.
  
  \begin{lemone} \label{Lemma:Limits in CSS}
  	\cite{rasekh2017left}
  	Let $G: I \to \W$ be a diagram. Then for an object $L$ in $\W$ the following are equivalent.
  	\begin{enumerate}
  		\item $L$ is the limit of $G$.
  		\item $\W_{/G}$ has a final object, which is a cone with cone point $L$.
  		\item  There exists an equivalence of complete Segal spaces $W_{/L} \simeq \W_{/G}$.
  	\end{enumerate}
   \end{lemone}

   There is one particular instance the previous lemma that we will need later on.
   
   \begin{corone} \label{Cor:Limit of arrow}
   	
   	Let $G: F(1) \to \W$ be a diagram. Then, the object $G(0)$ is the limit of $G$ and so there is an equivalence of complete Segal spaces 
   	$$\W_{/G(0)} \simeq \W_{/G}.$$
   \end{corone}

   \begin{defone} \label{Def:Arrow CSS}
   If $\W$ is a complete Segal space we denote the {\it complete Segal space of arrows} as $\W^{F(1)}$ (using \cref{Not:Free Arrow CSS}).
   It comes with a target map
   $$ t: \W^{F(1)} \to \W$$
   which is a Cartesian fibration if and only if $\W$ has finite limits. Concretely, this Cartesian fibration models the functor 
   \begin{center}
   	\begin{tikzcd}[row sep=0.01in]
   	 \W_{/-}: &[-0.4in]  \W^{op} \arrow[r]  & \CSS \\
   	 & x \arrow[r, mapsto] & \W_{/x} 	
   	\end{tikzcd}
   \end{center}
   where functoriality is given by pullback.
  \end{defone}
 
 We will require an important sub-fibration of $\W^{F(1)}$ that we will review here.
 
 \begin{defone} \label{Def:O CSS}
 	Let 
 	$$\O_\W \hookrightarrow \W^{F(1)}$$
 	be the subcategory with the same objects and morphisms pullback squares. Then, the composition $t:\O_\W\to \W$ is in fact a {\it right fibration} which corresponds to the functor 
 	  \begin{center}
 		\begin{tikzcd}[row sep=0.01in]
 			(\W_{/-})^{core}: &[-0.4in]  \W^{op} \arrow[r]  & \s \\
 			& x \arrow[r, mapsto] & (\W_{/x})^{core} 	
 		\end{tikzcd}
 	\end{center}
 \end{defone}

 \begin{remone}
 	The fibration $\O_\W$ has been covered extensively in the context of quasi-categories \cite[Section 6]{lurie2009htt}.
 \end{remone}

 \subsection{Grothendieck \texorpdfstring{$(\infty,1)$}{(oo,1)}-Topos}
 Grothendieck topos theory is very ubiquitous in algebraic geometry: Grothendieck $1$-toposes in classical algebraic geometry and 
 Grothendieck $(\infty,1)$-topos theory in derived algebraic geometry. However, we only focus on the fact they are special 
 cases of their elementary counterparts and thus we only give a minimal review.

 \begin{defone}[\citerefive{Proposition 2.2}] \label{Def: Groth Topos}
  A {\it Grothendieck topos} is a locally presentable $1$-category that satisfies weak descent.
 \end{defone}

 \begin{defone}[\citerefive{Theorem 6.9}] \label{Def:Groth Infty Topos}
  A {\it Grothendieck $(\infty,1)$-topos} is a presentable $(\infty,1)$-category that satisfies descent.
 \end{defone}
 
 \begin{exone} \label{Ex:S a Groth}
  The most simple example is $\s$, the $(\infty,1)$-category of Kan complexes. 
 \end{exone}
 
 We only need the following important observation relating Grothendieck $1$-toposes and $(\infty,1)$-toposes.
 
 \begin{corone}[\citerefive{Proposition 11.2}] \label{Cor Tau Groth Groth}
  Let $\G$ be a Grothendieck $(\infty,1)$-topos. Then the subcategory of $0$-truncated objects, denoted $\tau_0 \G$, 
  is a Grothendieck $1$-topos. 
 \end{corone}
 
 We will not require (or mention) any further details about Grothendieck topos theory.
 We refer the interested reader to \cite{maclanemoerdijk1994topos} for Grothendieck $1$-topos theory and  
 \cite{lurie2009htt} (using quasi-categories) or \cite{rezk2010toposes} (using model categories)
 for Grothendieck $(\infty,1)$-topos theory. 

 \subsection{Elementary \texorpdfstring{$(\infty,1)$}{(oo,1)}-Topos}
 We will assume familiarity with elementary $(\infty,1)$-topos theory later on and only review basic definitions.
 We give one detailed example (\cref{Ex:S EHT}) with the hope of giving the interested reader an intuition.
 The main source for elementary $(\infty,1)$-topos theory is \cite{rasekh2018elementarytopos}.
 
 \begin{defone}[\cite{lurie2009htt}] \label{Def:SOC}
 	Let $\W$ be a complete Segal space with finite limits. We denote by 
 	$$\Sub_\W: \W^{op} \to \set$$
 	the functor that takes every object $x$ to the set of (isomorphism classes) of subobjects $\Sub(x)$.
 	
 	We say $\Omega$ is a {\it subobject classifier} if it represents $\Sub$, meaning there is a natural bijection 
 	$$\Sub \cong \Map_\W(-,\Omega).$$ 
 \end{defone}

 \begin{defone}[\cite{rasekh2018elementarytopos}]
  Let $\W$ be a complete Segal space with finite limits. A map $p_{\U}: \U_* \to \U$ 
  is a {\it universe} if the induced map of right fibrations 
  $$\W_{/ \U} \hookrightarrow \O_\W$$
  is an embedding. Here $\O_\W$ denotes the right fibration as discussed in \cref{Def:O CSS}.
 \end{defone}

 \begin{defone}[\cite{rasekh2018elementarytopos}] \label{Def:Sufficient Universes}
  Let $\W$ be a complete Segal space with finite limits. We say $\W$ has {\it sufficient universes} if there exists a 
  collection of universes $\{ \U \}_{i \in I}$ such that the embeddings $\W_{/ \U} \to \O_\W$ are jointly surjective.
 \end{defone}

 \begin{remone}
  Intuitively this is telling us that for any morphism $f: A \to B$ in $\W$ there exists a universe $\U$ such that 
  $f$ is a pullback of $p_{\U} : \U_* \to \U$
  \begin{center}
   \pbsq{A}{\U_*}{B}{\U}{}{f}{p_{\U}}{}
  \end{center}.
 \end{remone}

 \begin{defone}[\cite{rasekh2018elementarytopos}] \label{Def EHT CSS}
  Let $\W$ be a complete Segal space. We say $\W$ is an {\it elementary complete Segal topos} if the following hold:
  \begin{enumerate}
   \item It has finite limits and colimits (\cref{Lemma:Limits in CSS}).
   \item It is locally Cartesian closed (\cref{Def:LCCC}).
   \item It has subobject classifier (\cref{Def:SOC}).
   \item It has sufficient universes (\cref{Def:Sufficient Universes}).
  \end{enumerate}
 \end{defone}

 The definition we gave only applies to complete Segal spaces, however we will need a definition for other models of $(\infty,1)$-categories as well.
 
 \begin{defone} \label{Def EHT QCat and sCat}
  A quasi-category $\Q$ is an {\it elementary quasi topos} if the complete Segal space $\Gamma(\Q)$ (or equivalently $t^!(\Q)$)
  is an elementary complete Segal topos.
  Moreover, a Kan enriched category $\K$ is an {\it elementary Kan enriched topos} if the quasi-category $N_\Delta(\K)$ is an  
  elementary quasi topos.
 \end{defone}

 \begin{remone}
  We say $\E$ is an elementary $(\infty,1)$-topos when we want to refer to the definition in any of those three models, without 
  specifying which model.
 \end{remone}
  
  We will need one result about elementary $(\infty,1)$-toposes later on: the {\it fundamental theorem of topos theory}.
  
  \begin{propone}[\cite{rasekh2018elementarytopos}] \label{Prop:Overcat EHT}
  	Let $\E$ be an elementary $(\infty,1)$-topos and $X$ an object. Then $\E_{/X}$ is an elementary $(\infty,1)$-topos with subobject classifier $\pi_2: \Omega\times X \to X$ and universes $\pi_2: \U \times X \to X$.
  \end{propone}

 \begin{remone}
 	The definition of elementary $(\infty,1)$-topos given here differs slightly from the definition in \cite{monaco2021eht} in the sense that there the universes are assumed to be closed under certain constructions (such as dependent products).
 \end{remone}

 This definition can be seen as a generalization of a Grothendieck $(\infty,1)$-topos as well as elementary topos. 
 
 \begin{defone}[\cite{maclanemoerdijk1994topos}]
  An {\it  elementary topos} is a locally Cartesian closed category with subobject classifier.
 \end{defone}

 \begin{propone}[\cite{rasekh2018elementarytopos}]
  Every Grothendieck $(\infty,1)$-topos is an elementary $(\infty,1)$-topos.
 \end{propone}
 
 \begin{propone}[\cite{rasekh2018elementarytopos}] \label{Prop:Underlying ET}
  Let $\E$ be an elementary $(\infty,1)$-topos. Then the subcategory of $0$-truncated objects, denoted $\tau_0\E$ is an elementary topos.
  We call it the underlying elementary topos.
 \end{propone}
 
 A general elementary $(\infty,1)$-topos does not have infinite colimits, as it is not an elementary condition. 
 We thus need an elementary alternative to infinite colimits that allow us to still recover some infinite constructions.
 This is achieved via the natural number object.
 
 \begin{defone}[\citenno{Definition 3.2}] \label{Def:NNO}
  A {\it natural number object} in an $(\infty,1)$-category $\E$ is an object $\mathbb{N}$ along with two morphisms 
  $s: \mathbb{N} \to \mathbb{N}$ and $o: 1 \to \mathbb{N}$, such that the triple $(\mathbb{N},s,o)$ is initial.
 \end{defone}

 Notice a natural number object is only useful in the context where we don't have infinite colimits, as the next example shows. 
  
 \begin{exone}[\citenno{Proposition 6.9}] \label{Ex:NNO Groth}
  If $\E$ is an elementary $(\infty,1)$-topos with countable colimits then the infinite coproduct $\coprod_{n \in \mathbb{N}} 1$ 
  is the natural number object. This in particular applies to every Grothendieck $(\infty,1)$-topos.
 \end{exone}

 This example gives us the following valuable recognition principle for elementary $(\infty,1)$-toposes that are not 
 Grothendieck $(\infty,1)$-toposes.
 
 \begin{corone} \label{Cor:NNO EHT Coprod}
  Let $\E$ be an elementary $(\infty,1)$-topos such that its natural number object $\mathbb{N}$ is not the infinite colimit $\coprod_{n \in \mathbb{N}} 1$.
  Then $\E$ is not a Grothendieck $(\infty,1)$-topos
 \end{corone}

 For more details on natural number objects for elementary $(\infty,1)$-topos see \cite{rasekh2021nno}.
 
 \par 
 
 In order to give a better understanding of the axioms of an elementary $(\infty,1)$-topos we give one detailed example.
 
 \begin{exone} \label{Ex:S EHT} 
  Let $\s_\bullet$ be the complete Segal space of Kan complexes 
  (we denote it $\s_\bullet$ rather than $\s$ in this example to make it clear it is a simplicial space).
  We already know that it is an elementary $(\infty,1)$-topos,
  as it is a Grothendieck $(\infty,1)$-topos (\cref{Ex:S a Groth}). 
  However, we want to use the fact that Kan complexes are well-understood to explain and give a better understanding of the 
  axioms of an elementary $(\infty,1)$-topos.
  \par 
  The existence of finite limits and colimits is a very standard condition and will not be discussed further.
   Similarly, it follows from our classical understanding of Kan complexes that $\s_\bullet$ is locally Cartesian closed. For further explanations see \cite{goerssjardine1999simplicialhomotopytheory, kapulkinlumsdaine2012kanunivalent}.
  
  Thus we move on to the existence of a subobject classifier. For that we first need to better understand mono maps in the complete 
  Segal space of Kan complexes. By definition a map of Kan complexes $f:X \to Y$ is mono if and only if the square 
  \begin{center}
   \pbsq{X}{X}{X}{Y}{id_X}{id_X}{f}{f}
  \end{center}
  is pullback square. This only holds if the map $f: X \to Y$ is a local equivalence i.e. a map $f$ that restricted to each path component 
  in $X$ becomes an equivalence. 
  \par 
  Equivalence classes of such maps are completely determined by a choice of path components in $Y$ (namely we choose the ones we want to be in $X$)
  and such a choice is determined by a map $Y \to \{ 0,1 \}$ (where we map the desired path components to $1$ and the rest to $0$).
  Thus $\{0,1\}$ is a subobject classifier in $\s$.
  \par 
  We now want to gain a better understanding of universes in $\s_\bullet$. Fix a large enough cardinal $\kappa$ and denote the full sub-complete Segal space of $\kappa$-small spaces by $\s^\kappa_\bullet$. Then we say a map of Kan 
  complexes is $\kappa$-small if each fiber is $\kappa$-small.
  The sub-space $\s^\kappa_0$ is itself an object in $\s_\bullet$. Thus for any Kan complex $K$ we can define the 
  space $\Map_\s(K, \s^\kappa_0)$. What are the points in this Kan complex?
  \par 
  A point in $\Map_\s(K, \s^\kappa_0)$ is a map of Kan complexes 
  $f:K \to \s^\kappa_0$. Every Kan complex $K$ is also a simplicial set and thus 
  $$K = \underset{\Delta^n \to K}{\colim} \Delta^n$$
  which means the map $f: K \to \s^\kappa_0$ is just a $K$-indexed family of maps $\Delta^n \to \s^\kappa_0 = \s^\kappa_{0n}$. 
  But $n$-cells in $\s^\kappa_0$ are just a choice of $n$ homotopic Kan complexes. In particular, a map $\Delta^0 \to \s^\kappa_0$
  is just a choice of $\kappa$-small Kan complex. 
  Thus a map $K \to \s^\kappa_0$ is a $K$-indexed diagram of $\kappa$-small Kan complexes, which we can denote 
  by a Kan fibration $\int_K f \twoheadrightarrow K$.
  \par 
  We can recover this map as the pullback of Kan complexes 
  \begin{center}
   \pbsq{\int_K f}{(\s_*)^\kappa_0}{K}{\s^\kappa_0}{}{}{}{}
  \end{center}
  where $(\s_*)^\kappa_0$ is the $0$-level of the complete Segal space of pointed $\kappa$-small Kan complexes, which comes with a forgetful map to 
  $\s^\kappa_0$.
  \par 
  This shows that $\s^\kappa_0$ is a universe in $\s_\bullet$. Using the fact that every map of Kan complexes is 
  $\kappa$-small for a large enough cardinal $\kappa$ proves that $\s_\bullet$ has sufficient universes.
 \end{exone}

 \begin{remone}
 	The example suggests that we should think of a universe as a ``small" copy of the core of the $(\infty,1)$-category sitting inside itself (the way that $\s^\kappa_0$ is sitting inside $\s_\bullet$). 
 \end{remone}

 \subsection{Filters}
 They key idea of the coming section is to use filters in the context of $(\infty,1)$-categories and $(\infty,1)$-toposes. 
 Hence, we will give a quick review of the history, applications and definitions we need.
 \par 
 Filters were first defined by Cartan in his study of topologies \cite{cartan1937filters, cartan1937ultrafiltres}. 
 However, it slowly found its way into many other branches of mathematics and in particular model theory, which is a development that happened in several stages. 
 \par 
 First, Skolem was trying to construct new models \cite{skolem1934earlyultraproducts, skolem1955models}, which led him to consider an early notion of ultra products. This was made further precise by \Los \ \cite{los1955ultraproduct}, who in particular proved the celebrated {\it \Los's theorem}, which is a powerful theorem helping us study ultra products of models. A formal definition of ultra products that is still in use came a little later \cite{fraynemorelscott1962ultraproduct}. For further discussions of ultra products see \cite{changkeisler1990modeltheory}. 
 \par 
 Given the unusual ability of the ultra product and \Los's theorem to relate finite concepts with infinite constructions it has found applications in many different and seemingly unrelated branches of mathematics. For examples, it was used to prove the Ax-Kochen isomorphism theorem \cite{axkochen1965diophantinei, axkochen1965diophantineii, axkochen1966diophantineiii}. However, it has also been used by Barthel, Schlank and Stapleton to further our understanding of chromatic homotopy theory \cite{bss2020chromatic}, using filters in the context of quasi-categories. 
 Our primary interest in filter quotients is from a topos theoretic perspective as has been discussed in \cref{Sec:Introduction}.
 \par
 We will only need the definition of a filter and so our review here will cover everything we need about a filter in the next sections. 
 
 \begin{defone} \label{Def:Filter}
  Let $(P, \leq)$ be a partially ordered set. A {\it filter} $F$ is a subset of $P$ that satisfies the following conditions.
  \begin{enumerate}
   \item $F \neq \emptyset$.
   \item $F$ is downward directed, meaning that for any two object $x, y \in F$ there exists $z \in F$ such that $z \leq x$ and $z \leq y$.
   \item $F$ is upward closed, meaning that if $x \leq y$ and $x \in F$, then $y \in F$.
  \end{enumerate}
 \end{defone}
 
 \begin{remone}
  Notice if $P$ has a maximum element, then every filter in $P$ necessarily includes that maximum element.
  Thus we could replace condition (1) with the condition that the maximum is in $F$.
 \end{remone}
  
  We have the following basic but crucial observation about filters that we will need in the next section.
  
  \begin{remone} \label{Rem:Filtered Colimits}
   If $F$ is a filter (and thus also a category), then the opposite category $F^{op}$ is a filtered category.
   This has two important implications about a diagram $D:F^{op} \to \C$ (where $\C$ is any cocomplete category) that we will use in the next section:
   \begin{enumerate}
   	\item The colimits of $D$ commutes with finite limits  \cite{maclane1998categories}.
   	\item Let $\F$ be a class of objects in $\C$ characterized by a right lifting property with respect to a set of monomorphisms with compact codomain.
   	If $D$ takes value in $\F$, then the colimit of $D$ is also in $\F$ 
   	For more details see the discussion of {\it small object arguments} in cofibrantly generated model categories \cite[2.1]{hovey1999modelcategories}\cite[11.1]{hirschhorn2003modelcategories}.
   	For an explicit argument in the context of simplicial sets see \cite[Theorem 4.1]{goerssjardine1999simplicialhomotopytheory}.
   \end{enumerate}
  \end{remone}

 \begin{defone} \label{Def:Ultrafilter}
  An {\it ultrafilter} $U$ of a poset $P$ is a filter of $P$ that is maximal, meaning there does not exist a filter $F$ such that 
  $U \not\subset F \not\subset P$.
 \end{defone}
 
 \begin{exone} \label{Ex:Principal Filter}
  Let $P$ be a poset and $x \in P$. Then the subset $\{y \in P: x \leq y \}$ is a filter. A filter of this form is called a {\it principal filter}.
 \end{exone}
 
 \begin{exone}
  Let $P$ be a poset with finite meets. Then any upward closed non-empty subset closed under finite meets is a filter.
 \end{exone}

 There is only one kind of filter we care about in the next section and so will discuss here in more detail.
 
 \begin{exone} \label{Ex:Filter Subob}
  Let $\E$ be a finitely complete $(\infty,1)$-category (for example an elementary $(\infty,1)$-topos).
  Then $\Sub(1)$ has the structure of a poset with finite meets. 
  Here the meet of two subobjects $U,V$ of $1$ is given by their product $U \times V \to 1$. 
  Given the explanation above, any subset of $\Sub(1)$ that includes $1$, is closed under finite products and is upward closed 
  is a filter in $\Sub(1)$. We will call any such filter a {\it product closed filter}.
 \end{exone}

 \begin{exone} \label{Ex:Filter on N}
  There is one interesting example of a filter we will use in \cref{Subsec:An Example}. 
  Let $\mathbb{N}$ be the set of natural numbers and $P(\mathbb{N})$ be the power set, which is a poset with the inclusion relation.
  Then cofinite sets (subsets with finite complement) form a product closed filter (commonly called the {\it Fr{\'e}chet filter}). 
  Indeed, if two subsets of $\mathbb{N}$ have 
  finite complement, then their intersection also has a finite complement. 
 \end{exone}
  
 \section{The Filter Quotient of \texorpdfstring{$(\infty,1)$}{(oo,1)}-Categories} \label{Sec:The Filter Quotient}
  The goal of this section is to construct a new elementary $(\infty,1)$-topos using a filter of $\Sub(1)$ (\cref{Ex:Filter Subob}).
 \par 
 In the first subsection we start more generally with a finitely complete $(\infty,1)$-category $\C$ 
 and use that to construct a new $(\infty,1)$-category $\C_\Phi$ along with a quotient functor from the original $(\infty,1)$-category 
 $P_\Phi : \C \to \C_\Phi$. We will give two main constructions depending on the model of 
 $(\infty,1)$-category:
 \begin{enumerate}
  \item One for Kan enriched categories, which leaves the collection of objects untouched and thus resembles the 
  definition of the filter quotient for elementary $1$-toposes \cite[Example A2.1.13]{johnstone2002elephanti} (\cref{Prop:E Phi simp Cat}).
  \item Another for complete Segal spaces and quasi-categories (\cref{Def:Filter Quotient CSS} \cref{The:Filter Quotient QCat}).
 \end{enumerate}
 We will then show these definitions agree with each other.
 \par 
 In the second section we prove that if $\E$ has finite limits, finite colimits, Cartesian closure, subobject classifiers, or universes, 
 then $P_\Phi$ will preserve them. This then implies that if $\E$ is an elementary $(\infty,1)$-topos then $\E_\Phi$ 
 is one as well.
 
 \begin{remone}
  Throughout this section the model of $(\infty,1)$-categories for $\E$ will change. The reader is advised to follow the 
  convention introduced in \cref{Subsec:Notation}.
  On other hand $\Phi$ will constantly denote a product closed filter in $\Sub(1)$.
 \end{remone}
 
 \subsection{The Filter Quotient Construction} \label{Subsec:The Filter Quotient}
 In this subsection we present two methods for constructing the filter quotient, 
 one that applies to Kan enriched categories and one for complete Segal spaces and quasi-categories.
 \par 
 First we will describe the filter quotient construction for Kan-enriched categories.
 Thus, let $\K$ be a finitely complete Kan enriched category. We will construct a new Kan enriched category which we denote by $\K_\Phi$.
 \par 
 The objects of $\K_\Phi$ are just the objects of $\K$. For the morphisms we first need some preliminary observations:
\begin{enumerate}
  \item $\Phi$ is a full subcategory of $\K$ and in particular there is a canonical functor 
 $\I: \Phi \hookrightarrow \K$.
 \item Let $- \times - : \K \times \K \to \K$ be the product functor. We can restrict it via $\I$ to a map 
 $$\I(-) \times - : \Phi \times \K \to \K$$
 \item We can apply opposite categories 
 $$ \I^{op}(-) \times - : \Phi^{op} \times \K^{op} \to \K^{op}$$
 \item We can product this functor with the identity functor
 $$ (\I^{op}(-) \times -) \times id_\K: \Phi^{op} \times \K^{op} \times \K \to \K^{op} \times \K$$
 \item Using the fact that $\K$ is a Kan-enriched category, we can post-compose the functor above
 $$\Map_\K(-,-) \circ ((\I^{op}(-) \times -) \times id_\K): \Phi^{op} \times \K^{op} \times \K \to \K^{op} \times \K \to 
 \Kan$$
 \item We can construct the adjoint of the functor above
 $$\F: \K^{op} \times \K \to Fun(\Phi^{op}, \Kan)$$
 defined as $\F(X,Y)(U) = \Map_\K(X \times U,Y)$.
 \item Now we observe that the diagram $\Phi^{op}$ is a filtered diagram in Kan complexes and so the colimit must be indeed a Kan complex, as Kan complexes are characterized via a right lifting property with respect to the maps $\Lambda^n_i \to \Delta^n$ (\cref{Rem:Filtered Colimits}).
 Thus we get a functor 
 $$\colim \circ \F: \K^{op} \times \K \to Fun(\Phi^{op},\Kan) \to \Kan$$
 \item We define 
 $$\Map_{\K_\Phi}(-,-) = \colim \circ \F$$
 and so in particular 
 $$\Map_{\K_\Phi}(X,Y) = \colim (\Map(- \times X,Y): \Phi^{op} \to \Kan)$$
\end{enumerate}
  
 First we prove that the Kan complexes $\Map_{\K_\Phi}(X,Y)$ give us indeed a new Kan enriched category.
 
 \begin{propone} \label{Prop:E Phi simp Cat}
  $\K_\Phi$ as defined above is a Kan enriched category.
 \end{propone}
 
 \begin{remone}
  Part of the proof here is a generalization of the argument given in \cite[Example A2.1.13]{johnstone2002elephanti}.
 \end{remone}

 \begin{proof}
  We use the observation that a simplicially enriched category is just a simplicial object in categories where each level has the same set 
  of objects. Thus it suffices to prove that for every natural number $k$, we get a category with object the objects $\K$ and morphisms sets
  $\Map_{\K_\Phi}(X,Y)_k$ and that composition respects simplicial boundary maps. 
  \par 
  First, recall that colimits in the category of simplicial sets are computed level-wise, thus we have 
  $$\Map_{\K_\Phi}(X,Y)_k = \colim (\Map(- \times X,Y)_k: \Phi^{op} \to \set)$$
  Let $X,Y,Z$ be three objects, we need a composition map 
  $$\Comp(X,Y,Z)_k: \Map_{\K_\Phi}(X,Y)_k \times \Map_{\K_\Phi}(Y,Z)_k \to \Map_{\K_\Phi}(X,Z)_k$$
  Let $[f]$ in $\Map_{\K_\Phi}(X,Y)_k$ and $[g]$ in $\Map_{\K_\Phi}(Y,Z)_k$ be two morphisms. 
  Then we can represent them as morphisms $f: X \times U \to Y$ and $g: Y \times V \to Z$.
  We now define the composition as
  $$[g] \circ [f] = [g(f \times id_V)]$$
  The colimit construction implies that the composition is well-defined. Moreover, as composition in $\K$ is associative this new composition 
  is associative as well. 
  \par 
  We have constructed a collection of categories which all have the same objects and where the hom set of the $k$-th category is given by $\Hom_{\K_\Phi}(X,Y)_k$.
  We need to show that this collection is a simplicial object in categories.
  However, this follows immediately from the fact that for every simplicial map $\delta: [k] \to [l]$ the 
  following diagram commutes: 
  \begin{center}
   \begin{tikzcd}[row sep=0.5in, column sep=1in]
    \Hom_{\K_\Phi}(X,Y)_l \times \Hom_{\K_\Phi}(Y,Z)_l \arrow[r, "\Comp(X \comma Y \comma Z)_l"] \arrow[d, "\delta_* \times \delta_*"] & 
    \Hom_{\K_\Phi}(X,Z)_l \arrow[d, "\delta_*"] \\
    \Hom_{\K_\Phi}(X,Y)_k \times \Hom_{\K_\Phi}(Y,Z)_k \arrow[r, "\Comp(X \comma Y \comma Z)_k"] & \Hom_{\K_\Phi}(X,Z)_k 
   \end{tikzcd}
  \end{center} 
   Thus we have proven that $\K_\Phi$ is a Kan enriched category.
  \end{proof}
  
  \begin{defone} \label{Def:Filter Quotient}
  Let $\K_\Phi$ be the category defined above. We call the Kan enriched category $\K_\Phi$ the {\it filter quotient of the Kan enriched category $\K$ 
  with filter $\Phi$}.
 \end{defone}

  Notice the category $\K_\Phi$ comes with a distinguished functor from $\K$.
  
  \begin{defone} \label{Def:P Phi Kan Enrich}
   There is a functor $P_\Phi: \K \to \K_\Phi$. It is the identity map on objects and it takes a morphism
   $f: X \to Y$ to the class $[f]$ of the morphism $f: X \times 1 \to Y$.
  \end{defone}
  
  We now will give an analogous construction for finitely complete complete Segal spaces (and by extension, quasi-categories). 
  The goal is to construct the filter quotient as a filtered colimit of complete Segal spaces. 
  This means we have to define our diagram and then prove the resulting colimit is a complete Segal space.
  Let $\W$ be a fixed complete Segal space with product closed filter $\Phi$ of $\Sub(1)$.
  
  \begin{defone}
   Let 
   $$\T_\Phi : \Phi^{op} \to \CSS$$
   be the diagram of complete Segal spaces that corresponds to the Cartesian fibration induced by the pullback
   \begin{center}
    \pbsq{\int \T_\Phi}{\W^{F(1)}}{\Phi}{\W}{}{}{}{}
   \end{center}
  \end{defone}
  
  \begin{remone}
   Concretely, the diagram $\Phi$ takes an object $V$ to the over-category $\W_{/V}$ and morphism $W \leq V$ to the product map 
   $$- \times W : \W_{/V} \to \W_{/W}$$
  \end{remone}
  
  We want to take the colimit of this diagram, but first we have to confirm that the colimit is a complete Segal space.
  
  \begin{lemone} \label{Lemma:Filtered Colimit CSS}
   Let $\Psi$ be a filtered category and let $F : \Psi \to \ssSet$ be a diagram in bisimplicial sets such that $F(X)$ is a 
   complete Segal space for each object $X$ in $\Psi$. Then the filtered colimit in $\ssSet$ is a complete Segal space. 
  \end{lemone}
  
  \begin{proof}
   As $\Psi$ is filtered, it suffices to observe that a complete Segal space is characterized via right lifting property against a class of monomorphisms with small codomain (\cref{Rem:Filtered Colimits}).
  
  A complete Segal space is by definition a {\it Reedy fibrant} simplicial space that satisfies the {\it Segal condition} and {\it completeness condition}. In \cite{rezk2001css, rezk2010thetanspaces} these conditions are explicitly described as right lifting properties with respect to these three sets of monomorphisms with small codomain:
  \begin{enumerate}
  	\item {\bf Reedy Fibrancy \cite[2.4]{rezk2001css}:}  The maps  
  	$$ F(n) \times \Lambda^l_k \underset{\partial F(n) \times \Lambda^l_k}{\coprod} \partial F(n) \times \Delta^l \to F(n) \times \Delta^l,$$
  	where $n \geq 0$, $l \geq 0$. Here $F(n)_{ij} = \Hom_{\Delta}([i],[n])$, $(\Delta^l)_{ij}=\Hom_{\Delta}([j],[l])$ and $\Lambda^l_k$ is the horn inside $\Delta^l$.
  	\item {\bf Segal Condition \cite[4.1]{rezk2001css}:}  The maps
  	$$ F(n) \times \partial \Delta^l \underset{G(n) \times \partial \Delta^l}{\coprod} G(n) \times \Delta^l \to F(n) \times \Delta^l,$$
  	where $n \geq 2$, $l \geq 0$. Here $G(n)$ is the {\it spine} inside $F(n)$.
  	\item {\bf Completeness Condition \cite[Proposition 10.1]{rezk2010thetanspaces}:}  The maps 
  	$$Z \times \partial \Delta^l \underset{\partial \Delta^l}{\coprod} \Delta^l \to Z \times \Delta^l,$$
  	where $l \geq 0$. Here $Z$ is defined as the pushout 
  	$$Z = F(3) \underset{F(1) \coprod F(1)}{\coprod} (F(0) \coprod F(0)).$$
  \end{enumerate}
  \end{proof}

  \begin{remone}
   An analogous (independent) proof for the special case of ultra products can be found in \cite[Lemma 3.13]{bss2020chromatic}.
  \end{remone}

  We now have all the necessary ingredients to define the filter quotient complete Segal space.
  
  \begin{defone} \label{Def:Filter Quotient CSS}
   Let $\W$ be a complete Segal space with finite limits. We define the {\it filter quotient} $\W_{\Phi}$ as 
   $$\W_{\Phi} = \underset{\Phi^{op}}{\colim} \ \T_\Phi$$
   and observe by the previous lemma that this is indeed a complete Segal space.
  \end{defone}
  
  The colimit construction also gives us the desired quotient functor.
  
  \begin{defone} \label{Def:P Phi CSS}
   Let $P_\Phi: \W \to \W_\Phi$ be the inclusion map into the filtered colimit, using the fact that $\W_{/1} = \W$ and $1 \in \Phi$.
  \end{defone}

  \begin{notone}
   In order to simplify notation we will usually denote the object $P_\Phi(X)$ in $\W_\Phi$ as
   $X$ again.
  \end{notone}

  We can give an analogous construction for quasi-categories. 
  
  \begin{theone} \label{The:Filter Quotient QCat}
   Let $\Q$ be a quasi-category and $\Phi$ be a product closed filter on $\Sub(1)$. 
   Then the simplicial set $\Q_{\Phi}$ defined as the colimit
   $$\Q_{\Phi} = \underset{\Phi^{op}}{\colim} \ \T_\Phi$$
   has the following properties:
   \begin{enumerate}
    \item It is a quasi-category.
    \item It is compatible with the definition for complete Segal spaces, in the sense that for a given complete Segal space $\W$
    we have equivalence 
    $$i^*_1(\W_\Phi) \simeq (i^*_1\W)_\Phi$$ 
   \end{enumerate}
   We call it the filter quotient of the quasi-category $\Q$ with respect to $\Phi$.
  \end{theone}
  
  \begin{proof}
   The key observation is that $\Gamma(\Q)$ is a complete Segal space with $\Gamma(\Q)_{n0} = \Q_n$ (\cref{Rem:Gamma of QCats}).
   Moreover, any filter $\Phi$ in $\Q$ is also a filter in $\Gamma(\Q)$. Thus, we can apply the filter 
   construction to get a complete Segal space $\Gamma(\Q)_\Phi$. 
   Restricting to the first row via $i_1^*(\Gamma(\Q))$ gives us a quasi-category. 
   But, in the construction of $\Gamma(\Q)_\Phi$ we constructed colimits level-wise, thus we have 
   $$i_1^*(\Gamma(\Q)_\Phi) = \underset{\Phi^{op}}{\colim} \  \T_\Phi  = \Q_{\Phi}$$
   This proves that $\Q_\Phi$ is a quasi-category and that it is compatible with the definition for complete Segal spaces.
  \end{proof}

  The question that remains is whether the construction for Kan enriched categories is compatible with the construction 
  for complete Segal spaces (the same way we just showed that the quasi-category and complete Segal space constructions are compatible).
  Concretely, let $\K$ be a Kan enriched category. Then can we compare $t^!N_\Delta(\K_{\Phi})$ 
  and $(t^!N_\Delta(\K))_\Phi$, noticing the fact that in the former we are using the Kan enriched filter quotient 
  construction, whereas in the latter we are using the complete Segal space filter quotient construction. 
  \par 
   Thus the final goal of this subsection is to prove that these constructions are indeed equivalent.
  For that we need an alternative characterization of the filter construction for Kan enriched categories.
  
  \begin{defone}
   Define an equivalence relation $\sim$ on the set of objects of $\K$ as follows
   $$X \sim Y \Leftrightarrow \text{there exists a } V \in \Phi \text{ such that } X \times V \simeq Y \times V$$
  \end{defone}

  \begin{defone} \label{Def:E Phi Quot}
   Let $Ob(\K^{quot}_\Phi)$ be a set of representatives of the equivalence classes of the equivalence relation given in the previous definition.
   Moreover, define $\K^{quot}_\Phi$ as the full subcategory of $\K_\Phi$ with set of objects in $Ob(\K^{quot}_\Phi)$.
  \end{defone}

  \begin{lemone} \label{Lemma:E Phi Quot equiv}
   The inclusion functor $\Ext: \K^{quot}_\Phi \hookrightarrow \K_\Phi$ is an equivalence.
  \end{lemone}

  \begin{proof}
   The functor is fully faithful by definition, thus it suffices to prove that the inclusion functor is essentially surjective.
   Let $X$ be an arbitrary object in $\K_\Phi$. Then there exists an object $Y$ in the set $Ob(\K^{quot}_\Phi)$ and $V \in \Phi$ 
   such that $X \times V \simeq Y \times V$. However, by definition of $\K_\Phi$, the two maps $\pi_1: X \times V \to X$ 
   and $id_X: X \to X$ are identified in $\Map_{\K_\Phi}(X,X)$, which implies that $X \simeq X \times V$.
   Similarly $Y \simeq Y \times V$. Thus we get 
   $$X \simeq X \times V \simeq Y \times V \simeq Y$$
   which proves that $\Ext$ is essentially surjective.
  \end{proof}
  
   \begin{lemone} \label{Lemma:Other Technical Lemma}
   Let $X,Y$ be two objects in $\K$. Then we have a bijection
   $$\Map_{\K_\Phi}(X,Y) \cong \underset{V \in \Phi^{op}}{\colim} \Map_\K(X \times V , Y \times V)$$
  \end{lemone}
  
  \begin{proof}
   It suffices to prove that for every $V$ in $\Phi$ we have a functorial bijection
   $$\Map_\K(X \times V , Y \times V) \xrightarrow{ \ \cong \ } \Map_{\K}(X \times V,Y)$$
   We have the following bijections:
   $$\Map_\K(X \times V , Y \times V) \cong \Map_\K(X \times V , Y) \times \Map_\K(X \times V , V) \cong \Map_\K(X \times V , Y)$$
   where the last bijection follows from the fact that $\Map_\K(X \times V , V)$ has either no elements or one element as $V$ is 
   a subobject of the final object and it is clearly not empty (we have $\pi_2: X \times V \to V$).
  \end{proof}

  We now want to give a construction of $\K_\Phi$ as a colimit analogous to our definition for complete Segal spaces. 
  The problem is that there is no straightforward way to construct over-categories for Kan enriched categories.
  Thus, we introduce a way to circumvent that problem. The key is the following observation from complete Segal spaces.
  
   \begin{remone}
   If $V$ is $(-1)$-truncated then the over-category $\W_{/V}$ is just a full subcategory of $\W$. Moreover the objects $X$ in $\W_{/V}$
   can be characterized by one of the following equivalent conditions:
   \begin{enumerate}
    \item There exists a map $X \to V$.
    \item The map  $\pi_1: X \times V \xrightarrow{ \ \simeq \ } X$ is an equivalence.
    \item There is an equivalence $X \simeq Y \times V$ for some object $Y$.
   \end{enumerate}
  \end{remone}
  
  This observation gives us the following definition:
  
  \begin{defone} \label{Def:ESV}
   Let $\K$ be a Kan enriched category and $V$ be a $(-1)$-truncated object. Then we denote by $\K^{S^V}$ the full 
   subcategory of $\K$ consisting of objects which satisfy $X \simeq X \times V$, where $X$ is an object in $\K$.
  \end{defone}

  \begin{defone}
   Let 
   $$\T_\Phi: \Phi^{op} \to \cat_\Delta$$
   be the Kan enriched functor defined as 
   $$\T_\Phi(V) = \K^{S^V}$$
  \end{defone}
  
  \begin{remone} \label{Rem:Filtered Colim CSS Map}
   In \cref{Def:CSS Mapping Space} we described the mapping space as a finite limit diagram. 
   Thus using the argument of \cref{Lemma:Filtered Colimit CSS} we see that for any filtered diagram $F:\Psi \to \ssSet$ 
   valued in complete Segal spaces we have an isomorphism of mapping Kan complexes
   $$\underset{V \in \Psi}{\colim} \Map_{F(V)}(X,Y) \cong \Map_{\underset{V \in \Psi}{\colim} F}(X,Y) $$
  \end{remone}
  
  We can now compare the filter construction for complete Segal spaces and Kan enriched categories.
  
  \begin{theone} \label{The:CSS and Kan En Agree}
   Let $(t^! N_\Delta) \circ \T_\Phi$ be the composition functor 
   $$(t^! N_\Delta) \circ \T_\Phi : \Phi^{op} \to \cat_\Delta \to \ssSet$$
   Then there is an equivalence of complete Segal spaces
   $$\underset{V \in \Phi^{op}}{\colim} (t^! N_\Delta) \circ \T_\Phi \xrightarrow{ \ \simeq \ } t^! N_\Delta(\K_\Phi).$$ 
  \end{theone}
  
  \begin{proof}
   We have already shown that there is an equivalence $\Ext: \K^{quot}_\Phi \to \K_\Phi$. 
   Thus it suffices to prove that there is an equivalence 
   $$\underset{V \in \Phi^{op}}{\colim} (t^! N_\Delta) \circ \T_\Phi \xrightarrow{ \ \simeq \ } t^! N_\Delta(\K^{quot}_\Phi)$$
   As a first step we want to define the map above by defining a cocone. 
   This means we have to construct functors $\K^{S_V} \to \K^{quot}_\Phi$ that are consistent with the functors $- \times V$.
   \par 
   Let 
   $$\F_V : \K^{S_V} \to \K^{quot}_\Phi$$
   be defined on objects by $\F_V(X \times V) = [X]$ and on mapping Kan complexes it is the inclusion map into the 
   colimit
   $$ \iota_V: \Map_{\K^{S_V}}(X \times V, Y \times V) \to \Map_{\K^{quot}_\Phi}([X],[Y]) = 
      \underset{W \in \Phi^{op}}{\colim} \Map_\K(X \times W, Y \times W) $$
   that takes a map $f: X \times V \to Y \times V$ to the class $[f]$ in the colimit.
   Here we used the alternative characterization of mapping Kan complexes in $\K_\Phi$ as proven in \cref{Lemma:Other Technical Lemma}.
   \par 
   We need to show that the collection of functors $\F_V$ give us a cocone over $\K^{quot}_\Phi$.
   However, this follows immediately from the fact that for any object $W \in \Phi$ with $W \leq V$, we have $[(X \times V) \times W] = [X \times W]$, 
   which proves that $\F_V(- \times W) = \F_W$. Thus the maps $\F_V$ give us a cocone. 
   \par 
   Applying the map $t^!N_\Delta$ we get a cocone of complete Segal spaces, 
   which gives us a universal map out of the colimit 
   $$\underset{V \in \Phi^{op}}{\colim} (\F_V): \underset{V \in \Phi^{op}}{\colim} t^! N_\Delta(\K^{S_V}) \rightarrow
   t^! N_\Delta(\K^{quot}_\Phi)$$
   We want to show that this map is an equivalence.
   \par 
   We know that $t^! N_\Delta(\K^{quot}_\Phi)$ and $t^! N_\Delta(\K^{S_V})$ are complete Segal spaces (\cref{Rem:NDelta sCat is CSS}). 
   Moreover, we prove in \cref{Lemma:Filtered Colimit CSS} that a filtered colimit of complete Segal spaces is a complete Segal space. 
   This proves that the left hand side is also a complete Segal space. 
   Thus it suffices to prove the map is Dwyer-Kan equivalence \cite[Theorem 7.7]{rezk2001css}. 
   \par 
   Clearly the map is a surjection on objects, thus we need to show we have an equivalence of mapping Kan complexes. 
   Fix two objects $X,Y$ in $\K$. Then we have a diagram of Kan complexes 
   \begin{center}
    \begin{tikzcd}[row sep=0.5in, column sep=0.5in]
     \underset{V \in \Phi^{op}}{\colim} \Map_{t^!N_\Delta(\K^{S_V})}(X \times V, Y \times V) \arrow[r] \arrow[d, "\simeq", "(1)"'] &
     \Map_{ t^!N_\Delta(\K^{quot}_\Phi)}([X],[Y]) \arrow[d, "\simeq", "(1)"'] \\
     \underset{V \in \Phi^{op}}{\colim} \Map_{\K^{S_V}}(X \times V, Y \times V) \arrow[r] \arrow[d, "\cong", "(2)"']&
     \Map_{\K^{quot}_\Phi}([X],[Y]) \arrow[d, "\simeq", "(3)"'] \\
     \underset{V \in \Phi^{op}}{\colim} \Map_\K(X \times V, Y \times V) \arrow[r, "\simeq", "(4)"'] &
     \Map_{\K_\Phi}(X,Y)    
    \end{tikzcd}
   \end{center}
   where the numbered morphisms are equivalences for the following reasons:
   \begin{enumerate}
    \item The map $t^!N_\Delta$ takes mapping Kan complexes to equivalent mapping Kan complexes (\cref{Rem:Mapping Spaces sCat to CSS}).
    \item $\K^{S_V}$ is a full subcategory $\K$, which gives us a bijection of mapping spaces (\cref{Def:ESV})
    \item The functor $\Ext$ is an equivalence of Kan enriched categories (\cref{Lemma:E Phi Quot equiv})
    \item This is an alternative characterization of mapping Kan complexes in $\K_\Phi$ (\cref{Lemma:Other Technical Lemma})
   \end{enumerate}
   This implies that the top horizontal map is also an equivalence of Kan complexes, which proves that
   we have a Dwyer-Kan equivalence of complete Segal spaces.
  \end{proof}

  \subsection{The Filter Quotient is an Elementary \texorpdfstring{$(\infty,1)$}{(oo,1)}-Topos} \label{Subsec:Filter Quotient EHT}
  The next step is to prove that if $\E$ is an elementary $(\infty,1)$-topos then $\E_\Phi$ is one as well.
  The key step of the proof is to show that $\E_\Phi$ also has sufficient universes.
  For that we need a better understanding of the filter quotient of the complete Segal space of cones.
  
  \begin{lemone} \label{Lemma:Cocones of E Phi}
   Let $I$ be a finite simplicial space (i.e. with finitely many non-degenerate cells) 
   and $\W$ a finitely complete complete Segal space with filter $\Phi$. 
   Let $\Phi$ be the induced filter on $\W^{I}$ consisting of constant functor $I \to \W$ with value $V \in \Phi$. 
   Then, we have 
   \begin{enumerate}
    \item an equivalence of functor complete Segal spaces 
    $$(\W^I)_{\Phi} \simeq (\W_{\Phi})^I$$
    \item an equivalence of cocones 
    $$(\W_{/I})_{\Phi} \simeq (\W_{\Phi})_{/I}$$
   \end{enumerate}   
  \end{lemone}
  
  \begin{proof}
   {\it (1)}
   First notice we have an equivalence 
   $$(\W^I)_{/V} \simeq (\W_{/V})^I$$
   as they both are the full subcategory of $\W^I$ consisting of diagrams that take value in the full subcategory $\W_{/V}$.
   Thus, we get the equivalence
   $$(\W^I)_\Phi = \underset{V \in \Phi^{op}}{\colim} (\W^I)_{/V} \simeq \underset{V \in \Phi^{op}}{\colim} (\W_{/V})^I 
   \simeq (\underset{V \in \Phi^{op}}{\colim} \W_{/V})^I = (\W_{\Phi})^I$$
   where the last equivalence follows from the fact that filtered colimits commute with exponents by finite simplicial spaces.
   \par 
   {(2)} We will now use our explicit description of cocones, the fact that filtered colimits commute with finite limits and the previous part to get the desired result:
   $$(\W_{/I})_{\Phi} \simeq (\W^{F(1) \times I} \times_{\W^{I}} * )_\Phi \simeq (\W^{F(1) \times I})_{\Phi} \times_{(\W^I)_\Phi} * \simeq 
   (\W_{\Phi})^{F(1) \times I} \times_{(\W_{\Phi})^I} *  = (\W_{\Phi})_{/I}$$
  \end{proof}
   
   Finally, we need a better understanding of $(\W_\Phi)_{/X}$ and $\O_{\W_\Phi}$.
   
  \begin{lemone}\label{Lemma:Over Phi}
  	Let $\W$ be a complete Segal space and let $\Phi$ be a filter of subobjects. Then for any object $X$, there is an equivalence
  	 $$(\W_{\Phi})_{/X} \simeq \underset{V \in \Phi^{op}}{\colim}(\W_{/X \times V}).$$
  \end{lemone}
  
  \begin{proof}
  	We have 
  	\begin{align*}
  		\underset{V \in \Phi^{op}}{\colim}(\W_{/X \times V}) & \simeq \underset{V \in \Phi^{op}}{\colim} (\W_{/V})_{/ \pi_2: X \times V \to V} & \text{\cref{Cor:Limit of arrow}}\\
  		& = \underset{V \in \Phi^{op}}{\colim} (\W_{/V})^{F(1)} \underset{\W_{/V}}{\times} F(0) & \text{\cref{Def:Over Cat CSS}} \\
  		&  \simeq (\underset{V \in \Phi^{op}}{\colim} \W_{/V})^{F(1)} \underset{\underset{V \in \Phi^{op}}{\colim} \W_{/V}}{\times} F(0) & \text{\cref{Lemma:Filtered Colimit CSS}/\cref{Lemma:Cocones of E Phi}} \\
  		& \simeq \W^{F(1)}_\Phi \underset{\W_\Phi}{\times} F(0)  = (\W_{\Phi})_{/X} & \text{\cref{Def:Filter Quotient CSS}}
  	\end{align*}
  \end{proof}

  \begin{lemone} \label{Lemma:O of Phi}
  	Let $\W$ be a complete Segal space and let $\Phi$ be a filter of subobjects. Then there is an equivalence of right fibrations 
  	$$\O_{\W_\Phi} \simeq \colim_{V \in \Phi^{op}} \O_{\W_{/V}}$$
  \end{lemone}

  \begin{proof}
  	By definition the right fibration $\O_{\W_\Phi}$ corresponds to the functor $((\W_\Phi)_{/-})^{core}$ and the right fibration $\colim_{V \in \Phi^{op}} \O_{\W_{/V}}$ corresponds to the functor $\colim_{V \in \Phi^{op}} (\W_{/V \times -})^{core}$, hence it suffices to prove the functors are naturally equivalent. Given the functoriality of colimit it suffices to prove that space $((\W_\Phi)_{/X})^{core}$ is equivalent to $\colim_{V \in \Phi^{op}}(\W_{/X \times V})^{core}$. However, this immediately follows from the previous lemma.
  \end{proof} 
	
  \begin{theone} \label{The E Phi EHT}
   Let $\W$ be a complete Segal space with finite limits. Then the functor 
   $$P_\Phi: \W \to \W_\Phi$$
   preserves 
   \begin{enumerate}
    \item finite limits and colimits
    \item subobject classifiers
    \item locally Cartesian closure
    \item universes
   \end{enumerate}
   In particular if $\E$ is an elementary $(\infty,1)$-topos, then $\E_\Phi$ is one as well.
  \end{theone}

  \begin{proof}
   We need to prove that $\W_\Phi$ satisfies the three conditions given in \cref{Def EHT CSS}. 
   We will confirm each separately.
   
   {\it Finite Limits and Colimits:}
   It suffices to prove that $\P_\Phi$ preserves finite limits and the case for finite colimits is analogous.
   We want to prove that the final object $1$ in $\W$ is also the final object in $\W_\Phi$. We have 
   $$\Map_{\W_\Phi}(X, 1) = \underset{V \in \Phi^{op}}{\colim} \Map_\W(X \times V, 1) \overset{(1)}{\simeq} 
   \underset{V \in \Phi^{op}}{\colim} \Delta[0] \overset{(2)}{=} \Delta[0]$$
   where are using the following facts:
   \begin{enumerate}
    \item $1$ is final in $\W$ and thus $\Map_\W(X \times U, 1)$ is contractible and the colimit construction is homotopy invariant.
    \item filtered colimits of the final object is again final.
   \end{enumerate}
   
   We now want to prove that $\W_\Phi$ has $I$-shaped limits. We have the following diagram 
   \begin{center}
    \begin{tikzcd}[row sep=0.5in, column sep=0.8in]
     \W_{/I} \arrow[r, "(P_{/I})_{\Phi}"] \arrow[dr, "(P_{\Phi})_{/I}"'] & (\W_{/I})_\Phi \arrow[d, "\simeq"] \\
     & (\W_{\Phi})_{/I}
    \end{tikzcd}
   \end{center}
   By the previous lemma the vertical map is an equivalence. By the previous paragraph $(P_{/I})_{\Phi}$ preserves final objects,  
   which then implies that $(P_{\Phi})_{/I}$ also preserves final objects.
   But a final object in cone category is just the limit and thus $P_\Phi$ preserves all finite limits.
   
   {\it Subobject Classifier:}
   Let $\Omega$ be the subobject classifier in $\W$. We want to prove that $\Omega$ is a subobject classifier in $\W_{\Phi}$.
   First, we have to determine the subobjects in $\W_\Phi$. Let $[f]$ be a morphism in $\Map_{\W_{\Phi}}(X,Y)$.
   Then $[f]$ is mono if and only if the following is a pullback square in $\W_\Phi$
   \begin{center}
    \pbsq{X}{X}{X}{Y}{id_X}{id_X}{\leb f \reb}{\leb f \reb}
   \end{center}
   However, we just proved that a pullback in $\W_\Phi$ is evaluated as a pullback of any representative. 
   Thus this is equivalent to the following being a pullback square 
   \begin{center}
    \comsq{X \times U}{X \times U}{X \times U}{Y}{id_{X \times U}}{id_{X \times U}}{f}{f}
   \end{center}
   where $f: X \times U \to Y$ is any representative of $[f]$.
   However, this is equivalent to $f: X \times U \to Y$ being mono in $\W$. 
   Thus we just proved that $[f]$ in $\W_\Phi$ is mono if and only if there exists a representative $f: X \times U \to Y$ that is mono in $\W$, 
   which we can state as a bijection 
   $$\Sub_{\W_\Phi}(X) \cong \underset{V \in \Phi^{op}}{\colim} \Sub_\W (X \times V)$$
   \par 
   We will use this to prove $\Omega$ is a subobject classifier. We have the following bijections
   $$\Map_{\W_\Phi}(X,\Omega) =  \underset{V \in \Phi^{op}}{\colim} \Map_\W(X \times V, \Omega \times V) \cong
   \underset{V \in \Phi^{op}}{\colim} \Map_\W(X \times V, \Omega) \cong
   \underset{V \in \Phi^{op}}{\colim} \Sub_\W (X \times V) \cong \Sub_{\W_\Phi}(X)$$
   where the first isomorphism follows from \cref{Lemma:Other Technical Lemma} and the last step follows from the bijection given in the last paragraph. 
   Thus $\Omega$ is also the subobject classifier in $\W_\Phi$.
   
   {\it Locally Cartesian Closure:}
   Let $f: X \times U \to Y \times U$ be a morphism in $\W_\Phi$. We already proved that $P_\Phi$ preserves finite limits, hence 
   it suffices to prove that if 
   $$f^*:\W_{/Y \times U} \to \W_{/X \times U}$$
   has a right adjoint, then the induced functor 
   $$f^*: (\W_\Phi)_{/Y \times U} \to (\W_\Phi)_{/X \times U}$$
   has a right adjoint.
   
   By \cite[Corollary 2.3]{gepnerkock2017univalence} it suffices to prove the following statement: For a morphism $g: Z \times U \to X \times U$ define the complete Segal space $(\W_\Phi)_{/f^*(-) \to Z \times U}$ as the pullback 
   \begin{center}
   	\begin{tikzcd}[row sep=0.5in, column sep=0.5in]
   		(\W_\Phi)_{/f^*(-) \to Z \times U} \arrow[r] \arrow[d] \arrow[dr, phantom, "\ulcorner", very near start] & (\W_\Phi)_{/Z \times U} \arrow[d, "g_*"] \\
   		(\W_\Phi)_{/Y \times U} \arrow[r, "f^*"] & (\W_\Phi)_{/X \times U} 
   	\end{tikzcd}.
   \end{center}
   Then $f^*$ has a right adjoint if and only if $(\W_\Phi)_{/f^*(-) \to Z \times U}$ has a final object. 
   
   Using the fact that filtered colimits commute with pullbacks, this pullback square is the filtered colimit of the pullback squares 
    \begin{center}
     \begin{tikzcd}[row sep=0.5in, column sep=0.5in]
   		\W_{/f^*(-) \to Z \times U} \arrow[r] \arrow[d] \arrow[dr, phantom, "\ulcorner", very near start] & (\W)_{/Z \times U} \arrow[d, "g_*"] \\
   		\W_{/Y \times U} \arrow[r, "f^*"] & (\W)_{/X \times U} 
   	 \end{tikzcd}.
    \end{center}
   By assumption $f^*: \W_{/Y \times U} \to \W_{/X \times U}$ has a right adjoint and so $\W_{/f^*(-) \to Z \times U}$ has a final object. 
   Hence, by the previous step, the filtered colimit, $(\W_\Phi)_{/f^*(-) \to Z \times U}$, also has a final object.
      
   {\it Universes:}
   Let $\U$ be a universe in $\W$. We want to prove that $\U$ is a universe in $\W_\Phi$, meaning there is an embedding
   $$(\W_\Phi)_{/\U} \hookrightarrow \O_{\W_\Phi}.$$
   First, recall that a filtered colimit of embeddings is an embedding (as being an embedding is a finite limit condition, which commutes with filtered colimits \cref{Rem:Filtered Colimits}). 
   Next, by \cref{Lemma:Over Phi}, we have 
   $$(\W_\Phi)_{/\U} \simeq \underset{V \in \Phi^{op}}{\colim} \W_{/V\times \U}$$
   and, by \cref{Lemma:O of Phi},
   $$\O_{\W_\Phi} \simeq \underset{V \in \Phi^{op}}{\colim} \O_{\W_{/V}}.$$
   Hence, it suffices to construct an embedding 
   $$\W_{/V\times \U} \hookrightarrow \O_{\W_{/V}}.$$
   However, the existence of this embedding follows from the fact that $\pi_2:\U \times V \to V$ is a universe in $\W_{/V}$ (\cref{Prop:Overcat EHT}). 
   \par 
   {\it $W_\Phi$ is an Elementary $(\infty,1)$-Topos:}
   Finally, we want to prove that if $\W$ is an elementary complete Segal topos, then $\W_\Phi$ is a complete Segal topos as well. 
   Following the previous steps, all that remains is to prove $\W_\Phi$ has sufficient universes. Let $[f]: X \to Y$ be an arbitrary map in $\W_\Phi$. 
   Then it can be represented by a map $f: X \times U \to Y \times U$ in $\W$. As $\W$ has sufficient universes $f$ is classified by 
   some universe $\U$. Given that pullbacks in $\W_\Phi$ are computed by pullbacks in $\W$, it follows that 
   $[f]$ is classified by $\U$ in $\W_\Phi$ as well. 
   \par 
   This proves that if $\E$ is an elementary complete Segal topos then $\E_\Phi$ is one as well. As we defined other models of 
   $(\infty,1)$-categories to be an elementary $(\infty,1)$-topos if the corresponding complete Segal space is one, the same 
   result holds for the other two models and hence we are done.
  \end{proof}

  \begin{remone}
   Some parts (such as preservation of limits and colimits and existence of adjoints) restricted to the case of ultra products were already proven in \cite[Lemma 3.17]{bss2020chromatic}.
  \end{remone}

   Notice there are trivial examples of filter quotients.
  
  \begin{exone} \label{Ex:Triv Filter}
   Let $\Phi$ be the minimal filter (which only includes $1$ itself). Then $\E_\Phi = \E$.
   On the other hand if $\Phi = \Sub(1)$, then $\E_\Phi$ is the trivial category with one object and identity map.
   This immediately follows from the fact that in this case the poset $\Phi^{op}$ has a final object and thus all colimits are
   just evaluation at that final object. However, that final object is just $0$ and $\E_{/0}$ is the trivial category.
  \end{exone}

  We finish this section by observing that this construction lifts the construction of the underlying toposes.
  
    \begin{theone} \label{The:E Phi Tau egal E Tau Phi}
   Let $\E$ be a locally Cartesian closed $(\infty,1)$-category with subobject classifier and 
   $\tau_0\E$ the elementary topos of $0$-truncated objects. 
   Moreover, let $\Phi$ be a product closed filter of $\Sub(1)$. Then there exists a functor 
   $$\tau_0(\E)_\Phi \to \tau_0(\E_\Phi)$$
   induced by $P_\Phi: \E \to \E_\Phi$, which is an equivalence of $1$-categories.
  \end{theone}

  \begin{proof}
   First we determine $0$-truncated objects in $\E_\Phi$. Let $X$ be an object in $\E_\Phi$. Then $X$ is $0$-truncated 
   if and only if the following commutative square 
   \begin{center}
    \pbsq{X}{X}{X}{X \times X}{id_X}{id_X}{\Delta_X}{\Delta_X}
   \end{center}
   is a pullback in $\E_\Phi$. This is equivalent to 
   \begin{center}
    \pbsq{X \times U}{X \times U}{X \times U}{(X \times X) \times U}{id_{X \times U}}{id_{X \times U}}{\Delta_X \times id_U}{\Delta_X \times id_U}
   \end{center}
   being a pullback in $\E$ for some object $U$ in $\Phi$.
   But the map $\Delta_X \times id_U$ is isomorphic to $\Delta_{X \times U}$ (as we have $U \cong U \times U$). 
   This means the square above being a pullback in $\E$ is equivalent to $X \times U$ being $0$-truncated in $\E$. 
   Thus we have proven that $X$ is $0$-truncated in $\E_\Phi$ if and only if it $X \times U$ is $0$-truncated in $\E$ for some $U$ in $\Phi$.
   \par 
   In particular, if $X$ is $0$-truncated in $\E$ then $P_\Phi(X)$ is an object in the subcategory $\tau_0(\E_\Phi)$, 
   which means we have the following commutative diagram
   \begin{center}
    \comsq{\tau_0(\E)}{\tau_0(\E_\Phi)}{\E}{\E_\Phi}{\tau_0(P_\Phi)}{}{}{P_\Phi}
   \end{center} 
   \par 
   Notice $\Phi$ is also a filter in $\tau_0(\E)$ and so we can apply the filter quotient construction to get a diagram
   \begin{center}
    \begin{tikzcd}[row sep=0.5in, column sep=0.5in]
     \tau_0(\E) \arrow[r, "\tau_0(P_\Phi)"] \arrow[d, "P_\Phi"] & \tau_0(\E_\Phi) \\
     \tau_0(\E)_\Phi \arrow[ur, dashed]
    \end{tikzcd}
   \end{center}
   We now want to show that we can lift the diagram to a functor $P_\Phi: \tau_0(\E)_\Phi \to \tau_0(\E_\Phi)$.
   \par 
   The map is already defined for objects, but we need to define it for $\Hom$ sets, which means we need a map 
   $$\Hom_{(\tau_0\E)_\Phi}(X,Y) \to \Hom_{\tau_0(\E_\Phi)}(X,Y)$$
   Using their definition as colimits this means we need a map 
   $$\underset{V \in \Phi^{op}}{\colim} \Hom_{\tau_0\E}(X \times V,Y) \to \underset{V \in \Phi^{op}}{\colim} \Hom_{\E}(X \times V,Y)$$
   But $\tau_0\E$ is a full subcategory of $\E$ and so we can take this map to be the identity.
   \par 
   We now want to prove that this induced functor is an equivalence of categories. 
   The functor is fully faithful by construction, so we only need to prove it is essentially surjective.
   Let $X$ be an object in $\tau_0(\E_\Phi)$. Then there exists a $U$ such that $X \times U$ is $0$-truncated.
   This means that $X \times U$ is an object in $\tau_0(\E)$ and thus also an object in $(\tau_0(\E))_\Phi$.
   Finally by construction of the filter quotient, the map $U \times X \to X$ is an equivalence. 
   Thus we have proven that every object $X$ is equivalent to an object in the image, namely $X \times U$.
  \end{proof}
  
  One interesting implication of this result is the preservation of natural number objects.
  
  \begin{theone} \label{The:P Phi preserves NNO}
   Let $\E$ be an elementary $(\infty,1)$-topos and $\mathbb{N}$ be the natural number object. 
   Moreover, let $\Phi$ be a product closed filter on $\E$. Then $P_\Phi(\mathbb{N})$ is the natural number object in $\E_\Phi$.
  \end{theone}

  \begin{proof}
   A natural number object $\mathbb{N}$ is always $0$-truncated and thus lives in $\tau_0(\E)$. 
   By \cite{johnstone2002elephanti} $P_\Phi(\mathbb{N})$ is thus also a natural number object in the filter quotient elementary topos $\tau_0(\E)_\Phi$.
   By the equivalence above, it is thus also a natural number object in $\tau_0(\E_\Phi)$.
   But by \cite[Lemma 3.13, Theorem 5.3]{rasekh2021nno} every natural number object in the underlying elementary topos is a natural number object in 
   the whole elementary $(\infty,1)$-topos. Thus $\E_\Phi$ has a natural number object, namely $P_\Phi(\mathbb{N})$.
  \end{proof}

  We can now combine this theorem with \cref{Cor Tau Groth Groth} to get a powerful method to construct
  many elementary $(\infty,1)$-toposes that are not Grothendieck $(\infty,1)$-toposes.
  
  \begin{corone} \label{Cor EHT not GHT}
   Let $\E$ be an elementary $(\infty,1)$-topos and $\Phi$ be a product closed filter such that the filter quotient elementary topos 
   $\tau_0(\E)_\Phi$ is not a Grothendieck topos. Then the filter quotient $\E_\Phi$ is an elementary $(\infty,1)$-topos that is not 
   a Grothendieck $(\infty,1)$-topos.
  \end{corone}
  
  \begin{remone}
   The conditions given in the previous corollary are sufficient conditions for $\E$ to be elementary without being Grothendieck. 
   It is not known whether the conditions are also necessary i.e. if there is an elementary $(\infty,1)$-topos that is not a
   Grothendieck $(\infty,1)$-topos, but for which the underlying elementary topos is a Grothendieck topos.
  \end{remone}
 
  In the next section we will focus on specific filter quotients and then use it, combined with the corollary above, to construct 
  elementary $(\infty,1)$-toposes are not Grothendieck $(\infty,1)$-toposes

\section{Filter Products} \label{Sec:Filter Products}
 In this section we want to restrict our attention to specific filter quotients: {\it  filter products}. 
 In the first subsection we study general filter product. In the second subsection 
 we focus on one specific example of a filter product.
 
 \subsection{\Los's Theorem for Equivalences} \label{Subsec:Los Theorem for Equivalences}
 Throughout this subsection $\C$ is a fixed finitely complete $(\infty,1)$-category such that $1$ only has two subobjects, $I$ is a set 
 and $\Phi$ is a product closed filter of $P(I)$, the power set of $I$.
 Notice $\C^I$ is also finitely complete.  We want to use the filter $\Phi$ of $P(I)$ to build a filter on $\C^I$.
 \par 
 First, we observe by assumption there is a bijection
 $$\{0,1 \} \cong \Sub_\C(1) $$
 Moreover, we have 
 $$\Sub_{\C^I}(1) = \Sub_\C(1)^I = \{0,1 \}^I = P(I)$$
 and so a filter $\Phi$ on $P(I)$ is automatically a filter on $\C^I$ and we can use it to define the filter quotient $(\C^I)_{\Phi}$.

 \begin{defone}
  Let $\C$ be a finitely complete $(\infty,1)$-category such that $1$ has only two subobjects, $I$ a set and $\Phi$ a product closed
  filter of $P(I)$. Then the resulting filter quotient $(\C^I)_\Phi$ is called the {\it filter product} and denoted 
  $\prod_\Phi \C$
 \end{defone}

 \begin{remone}
  If the filter $\Phi$ is an ultrafilter, then $\prod_\Phi \C$ is often called the {\it ultra product}.
 \end{remone}
 
  Our goal is to prove analogues of \Los's theorem \cite{los1955ultraproduct} for homotopies and equivalences in a filter product. 
  We will actually prove a more general statement about filter quotients, from which the desired results for 
  filter products will follow as an immediate corollary.
 
 \begin{lemone} \label{Lemma:Equal in Fil Colimits of Spaces}
  Let $\Phi^{op}$ be a filtered diagram and $F: \Phi^{op} \to \s$ be a diagram of spaces. 
  Then
  \begin{enumerate}
   \item We have isomorphism of sets 
   $$\pi_0(\colim F) \cong \colim ( \pi_0(F))$$
   \item Two points $x,y \in \colim(F)$ are homotopic if and only if there exists $S \in \Phi^{op}$ and $x', y' \in F(S)$ 
   such that $x' \simeq y'$ and $\iota(x') = x$, $\iota(y') = y$, where $\iota: F(S) \to \colim(F)$ is the universal cocone map.
  \end{enumerate}
 \end{lemone}

 \begin{proof}
  The first part follows immediately from the fact that $\pi_0$ commutes with filtered colimits. 
  Given that the second part focuses on the existence of equivalences we can restrict our attention to $\pi_0$ of that diagram.
  \par 
  Thus we have to prove that $x = y$ in $\colim (\pi_0(F))$ if and only if there exists $S$ and $x',y' \in \pi_0(F(S))$ 
  such that $\iota(x') = x, \iota(y') = y$ and $ x'= y'$. 
  However, that is just the definition of the colimit.
 \end{proof}
 
 \begin{theone} \label{The:Los Two maps Equiv}
  Let $\C$ be a finitely complete $(\infty,1)$-category and $\Phi$ a filter. Then for two morphisms $f, g: X \to Y$ in $\E_\Phi$ 
  $$ f \simeq g \text{ in } \C_\Phi \Leftrightarrow \exists U \in \Phi ( f \times id_U \simeq g \times id_U \text{ in } \C)$$
 \end{theone}
 
 \begin{proof}
  Observe that 
  $$\Map_{\C_\Phi}(X,Y) = \underset{V \in \Phi^{op}}{\colim} \Map_\C(X \times V, Y \times V)$$
  So, the result follows immediately from \cref{Lemma:Equal in Fil Colimits of Spaces} as two morphisms are equivalent 
  if and only if they are in the same path-component.
 \end{proof}

 Using the result for filter products we get the following corollary.
 
 \begin{corone} \label{Cor:Los Two maps Equiv}
  Let $(f_i)_{i \in I}, (g_i)_{i \in I} : (X_i)_{i \in I} \to (Y_i)_{i \in I}$ be two maps. 
  Then $(f_i)_{i \in I} \simeq (g_i)_{i \in I}$ in $\prod_\Phi \C$ if and only if 
  $$\{ i \in I : f_i \simeq g_i \text{ in } \C \} \in \Phi$$
 \end{corone}
 
 \begin{proof}
  Assume $\{ i \in I : f_i \simeq g_i \text{ in } \C \} \in \Phi$. Then we can use \cref{The:Los Two maps Equiv} with 
  $U = \{ i \in I : f_i \simeq g_i \text{ in } \C \}$ to deduce that  $(f_i)_{i \in I} \simeq (g_i)_{i \in I}$.
  
  On the other side, assume that $(f_i)_{i \in I} \simeq (g_i)_{i \in I}$. Then, by \cref{The:Los Two maps Equiv},
  there exists a $U \in \Phi$ such that for all $i \in U$, $f_i \simeq g_i$, meaning that 
  $$U \subseteq \{ i \in I : f_i \simeq g_i \text{ in } \C \}.$$
  However, by definition of a filter (\cref{Def:Filter}), this implies that $\{ i \in I : f_i \simeq g_i \text{ in } \C \} \in \Phi$.
 \end{proof}

 \begin{theone} \label{The:Los Maps Equiv}
  Let $\C$ be a finitely complete $(\infty,1)$-category, $\Phi$ a filter on $\C$. 
  Then a map $f$ in $\C_\Phi$ is an equivalence if and only if there exists $U \in \Phi$ such that $f \times id_U$ is an equivalence in $\C$.
 \end{theone}
 
 \begin{proof}
  Let $\C^{core}$ be the underlying $(\infty,1)$-groupoid of $\C$. Then we have an equivalence 
  $$(\C_\Phi)^{core} \simeq  (\C^{core})_\Phi$$
  Indeed, this follows from the fact that in the complete Segal space model, we have $\C^{core} = \C_0$ (\cref{Def:Core of CSS})
  and the filter quotient is defined as a level-wise colimit (\cref{Def:Filter Quotient CSS}). 
  \par 
  A morphism $f$ in $\C_\Phi$ is an equivalence if and only if it is in $(\C_\Phi)^{core}$, which is equivalent to $(\C^{core})_\Phi$.
  But, by definition of filtered colimits, this is equivalent to $f \times id_U$ being invertible for some $U \in \Phi$.
 \end{proof}
 
 Again, we can restrict our attention to filter products.

 \begin{corone} \label{Cor:Los Maps Equiv}
  Let $\C$ be a finitely complete $(\infty,1)$-category such that $1$ has two subobjects, $I$ a set and $\Phi$ a filter on $P(I)$. 
  Then a map $(f_i)_{i \in S}$ in $\prod_\Phi \C$ is an equivalence if and only if 
  $$ \{ i \in I : f_i \text{ is an equivalence } \} \in \Phi$$
 \end{corone}
 
 \begin{proof}
  Assume $\{ i \in I : f_i \text{ is an equivalence } \} \in \Phi$. Then we can use \cref{The:Los Maps Equiv} with 
  $U = \{ i \in I : f_i \text{ is an equivalence } \}$ to deduce that  $(f_i)_{i \in I}$ is an equivalence.
 
  On the other side, assume that $(f_i)_{i \in I}$ is an equivalence. Then, by \cref{The:Los Maps Equiv},
  there exists a $U \in \Phi$ such that for all $i \in U$, $f_i$ is an equivalence, meaning that 
  $$U \subseteq \{ i \in I : f_i \text{ is an equivalence } \}.$$
 However, by definition of a filter (\cref{Def:Filter}), this implies that $\{ i \in I : f_i \text{ is an equivalence } \} \in \Phi$.
 \end{proof}

 \begin{remone}
  We can use these results to characterize truncated objects in a filter quotient. 
  See \cite[Subsection 6.2]{rasekh2018truncations} for more details.
 \end{remone}

 \begin{remone}
  We only added the restriction that $1$ in $\C$ has only two subobject in order to guarantee that $\Phi$ is also a filter on $\C^I$. 
  With enough care that condition could possibly be relaxed.
 \end{remone}

 \begin{remone}
  Notice in the actual \Los's theorem, the filter needs to be an ultrafilter, which we did not assume.
  That is because \Los's theorem (as used in model theory) holds for all formulas. In particular, it holds for formulas that include
  negations, which require ultrafilters, as they have a certain closure property under set complements. 
  \par 
  On the other hand we prove very particular results, none of which include any negation and that is why these statements hold 
  for any filter. I want to thank Peter Lumsdaine for making me aware of this fact.
 \end{remone}

 We have now gathered enough background to give examples of non-presentable $(\infty,1)$-toposes.
   
\subsection{Examples of Filter-Quotients that are not Grothendieck \texorpdfstring{$(\infty,1)$}{(oo,1)}-Toposes} \label{Subsec:An Example}
 In this subsection we finally give examples of elementary $(\infty,1)$-toposes that are not Grothendieck toposes. 
 First we use our general knowledge to construct a whole class of examples. 
 Then we focus on one specific example and make the construction as explicit as possible.
 \begin{exone} \label{Ex:KanS Phi not GHT}
   As before, let $\Kan$ be the Kan enriched category of Kan complexes, let $I$ be a set and $\Phi$ a filter on $P(I)$.
   Notice that we indeed have $\Sub(1) = \{ 0 , 1 \}$. Thus, we can define the filter product $\prod_\Phi \Kan$
   and by \cref{The E Phi EHT} it is still an elementary $(\infty,1)$-topos. By 
   \cref{The:E Phi Tau egal E Tau Phi} the underlying elementary topos is $\prod_\Phi \set$.
   Based on \cref{Cor EHT not GHT} we only need to show that $\prod_\Phi \set$ is not a Grothendieck topos. 
   \par 
   For example let $I$ be a set and $\Phi$ be a non-principal filter. Then $\prod_\Phi \set$ is not a 
   Grothendieck topos \cite[Theorem 3.4]{adelmanjohnstone1982serreclasses}. Thus there are at least as many elementary $(\infty,1)$-toposes that are not 
   Grothendieck $(\infty,1)$-toposes as there are non-principal filters of sets. 
  \end{exone}
  
 Here is one particular example that satisfies the conditions given in \cref{Ex:KanS Phi not GHT}: 
 Let $\Phi$ be the filter of cofinite subsets of $\mathbb{N}$, the set of natural numbers (\cref{Ex:Filter on N}). We can thus apply the 
 filter quotient to the topos $\Kan^\mathbb{N}$. The resulting topos $\prod_\Phi \Kan$ is not a Grothendieck $(\infty,1)$-topos 
 as its subcategory of $0$-truncated objects is not a Grothendieck topos \cite[Example D5.1.7]{johnstone2002elephanti}.
 \par 
 We end this section by giving a detailed description of this particular example 
 without using the language of filter quotients. 
 Thus a reader who only wants to see an example can avoid the technical details of the previous section.
 We will refrain from giving detailed proofs here and refer the interested reader to the proofs in the previous section.
 \par 
 Let $\Kan$ be the Kan enriched category of Kan complexes. Then $\Kan^\mathbb{N}$ is the Kan enriched category with
 \begin{enumerate}
  \item objects tuples $(X_n)_{n \in \mathbb{N}}$, where $X_n$ is a Kan complex, and 
  \item morphisms level-wise morphisms
  $$\Map_{\Kan^\mathbb{N}}((X_n)_{n \in \mathbb{N}},(Y_n)_{n \in \mathbb{N}}) = \prod_{n \in \mathbb{N}} \Map_{\Kan}(X_n,Y_n)$$
 \end{enumerate}
 We will now define a new Kan enriched category $\prod_\Phi \Kan$. However for that we need an equivalence relation
 on the mapping Kan complexes $\Map_{\Kan}(X_n,Y_n)$. We will define an equivalence relation for each level of the Kan complex 
 and then show that the simplicial operator maps respect the equivalence relation, thus giving us a new simplicial set.
 \par 
 Denote the set of natural numbers bigger than $m$ by $\mathbb{N}_{\geq m}$.
 We define an equivalence relation on the set
 $\coprod_{m \in \mathbb{N}} (\Map_{\Kan^\mathbb{N}}((X_n)_{n \in \mathbb{N}_{\geq m}},(Y_n)_{n \in \mathbb{N}_{\geq m}}))_k$ as follows.
 Let 
 $$(f_n)_{n \in \mathbb{N}_{\geq m_1}} \sim (g_n)_{n \in \mathbb{N}_{\geq m_2}} \Leftrightarrow 
 \text{ there exists } N > m_1,m_2 \text{ such that for all } n > N: (f_n = g_n)$$
 Notice if we have two maps 
 $$f,g: \Delta[k] \times X_n \to Y_n$$
 for $n > m_1$ and $f \sim g$ then $f\delta \sim g\delta$ for any simplicial map $\delta : \Delta[l] \to \Delta[k]$, 
 as $f_n\delta = g_n\delta$ for $n$ large enough. Thus imposing the equivalence relation level-wise still gives us a simplicial set, 
 and in fact a Kan complex.
 \par 
 We can use that to define a new category $\prod_\Phi \Kan$ with
 \begin{enumerate}
  \item objects tuples $(X_n)_{n \in \mathbb{N}}$ (so the same objects as $\Kan^\mathbb{N}$) and
  \item morphisms level-wise equivalence classes of morphisms
  $$\Map_{\prod_\Phi \Kan}((X_n)_{n \in \mathbb{N}},(Y_n)_{n \in \mathbb{N}}) = 
  \left[ \coprod_{m \in \mathbb{N}} \prod_{n \in \mathbb{N}_{\geq m}} \Map_{\Kan}(X_n,Y_n) \right] / \sim .$$
 \end{enumerate}
  Intuitively a $k$-cell in $\Map_{\prod_\Phi \Kan}((X_n)_{n \in \mathbb{N}},(Y_n)_{n \in \mathbb{N}})$
  is a class 
  $$[f_n]: (X_n)_{n \in \mathbb{N}_{\geq m}} \times \Delta[k] \to (Y_n)_{n \in \mathbb{N}_{\geq m}}$$
  where two morphisms $f_n$, $g_n$ are in the same class if $f_n= g_n$ for $n$ large enough.
  \par 
  We need to confirm that we actually get a category. It suffices to check that we have a composition. 
  For two morphisms classes $[f_n]: (X_n)_{n \in \mathbb{N}_{\geq m_1}} \to (Y_n)_{n \in \mathbb{N}_{\geq m_1}}$ and 
  $[g_n]: (Y_n)_{n \in \mathbb{N}_{\geq m_2}} \to (Z_n)_{n \in \mathbb{N}_{\geq m_2}}$ we define the composition as
  $$[g_n \circ f_n]: (X_n)_{n \in \mathbb{N}_{\geq Max(m_1,m_2)}} \to (Z_n)_{n \in \mathbb{N}_{\geq Max(m_1,m_2)}}.$$
  and notice this definition of composition is indeed well-defined.
  \par 
  We now want to prove that $\prod_\Phi \Kan$ is an elementary $(\infty,1)$-topos. 
  We need to show that it has finite limits and colimits, is locally Cartesian closed, has subobject classifier and universes.
  We will treat each separately.
  \par 
  First, $\prod_\Phi \Kan$ has a final object, namely $(\Delta[0])_{n \in \mathbb{N}}$.
  Now, for two morphisms $(f_n)_{n \in \mathbb{N}_{\geq m_1}}: (X_n)_{n \in \mathbb{N}_{\geq m_1}} \to (Z_n)_{n \in \mathbb{N}_{\geq m_1}}$ 
  and $(g_n)_{n \in \mathbb{N}_{\geq m_2}}: (Y_n)_{n \in \mathbb{N}_{\geq m_2}} \to (Z_n)_{n \in \mathbb{N}_{\geq m_2}}$
  a routine calculation shows that the pullback is just the level-wise pullback 
  $(X_n \times_{Z_n} Y_n)_{n \in \mathbb{N}_{\geq Max(m_1,m_2)}}$.
  This proves that $\prod_\Phi \Kan$ has finite limits. The argument for finite colimits is analogous.
  \par 
  Next, observe that $\prod_\Phi \Kan$ is locally Cartesian closed. Again we use the Cartesian closure of Kan complexes level-wise. 
  Concretely, for a map $(p_n)_{n \in \mathbb{N}_{\geq m_1}}:(Z_n)_{n \in \mathbb{N}_{\geq m_1}} \to (X_n)_{n \in \mathbb{N}_{\geq m_1}}$ and 
  $(f_n)_{n \in \mathbb{N}_{\geq m_2}}:(X_n)_{n \in \mathbb{N}_{\geq m_2}} \to (Y_n)_{n \in \mathbb{N}_{\geq m_2}}$, the map $f_*p$ is given by
  $$((f_n)_*(p_n))_{n \in \mathbb{N}_{\geq Max(m_1,m_2)}} : ((f_n)_*(Z_n))_{n \in \mathbb{N}_{\geq Max(m_1,m_2)}} \to (Y_n)_{n \in \mathbb{N}_{\geq Max(m_1,m_2)}}$$
  which gives us the desired right adjoint to the pullback.
  \par 
  Next we want to find the subobject classifier. We showed in \cref{Ex:S EHT} that the subobject classifier in $\Kan$ is 
  $\{0,1\}$. We will now show that the constant sequence $(\{0,1\})_{n \in \mathbb{N}}$ is the subobject classifier in $\prod_\Phi \Kan$.
  First, notice that a map 
  $(f_n)_{n \in \mathbb{N}_{\geq m}}: (X_n)_{n \in \mathbb{N}_{\geq m}} \to (Y_n)_{n \in \mathbb{N}_{\geq m}}$
  is mono if and only if $f_n$ is mono for $n > N$ for some $N$.
  But every mono $f_n: X_n \to Y_n$ is uniquely determined by a map $Y_n \to \{0,1\}$ as it is a subobject classifier in $\Kan$,
  This means we get maps $Y_n \to \{ 0,1 \}$ for $n > N$, which is exactly the data of a map 
  $(Y_n)_{n \in \mathbb{N}} \to (\{0,1\})_{n \in \mathbb{N}}$. 
  So, every mono map in $\prod_\Phi \Kan$ is determined by a map into $(\{0,1\})_{n \in \mathbb{N}}$ and we can 
  show this assignment is unique, proving it is a subobject classifier.
  \par 
  Finally, we need to show $\prod_\Phi \Kan$ has universes. 
  Recall that the universes in $\Kan$ are the objects $\s^\kappa_0$, where $\kappa$ is a cardinal 
  (see \cref{Ex:S EHT} for more details).
  We want to show that the objects $(\s^\kappa_0)_{n \in \mathbb{N}}$ give us universes in $\prod_\Phi \Kan$.
  Let 
  $(f_n)_{n \in \mathbb{N}_{\geq m}}: (X_n)_{n \in \mathbb{N}_{\geq m}} \to (Y_n)_{n \in \mathbb{N}_{\geq m}}$ be an arbitrary map.
  Choose a cardinal $\kappa$, such that every morphism $f_n$ is $\kappa$-small.
  Then every map $f_n$ is a pullback of a map $Y_n \to \s^\kappa$. However, we just showed that pullbacks 
  in the category $\prod_\Phi \Kan$ are evaluated level-wise.
  Thus $(f_n)_{n \in \mathbb{N}_{\geq m}}$ is the pullback of $(Y_n)_{n \in \mathbb{N}_{\geq m}} \to (\s^\kappa)_{n \in \mathbb{N}_{\geq m}}$,
  which shows that every map is classified by a universe. We can show that this assignment gives us an equivalence. 
  \par 
  In order to finish this example we need to show that $\prod_\Phi \Kan$ is not a Grothendieck $(\infty,1)$-topos. 
  Following \cref{Cor:NNO EHT Coprod} it suffices to prove that its natural number object is not equivalent to the countable colimit.
  \par 
  The argument we give here is analogous to \cite[Example D5.1.7]{johnstone2002elephanti}.
  The natural number object in $\prod_\Phi \Kan$ is the constant sequence $(\mathbb{N})_{n \in \mathbb{N}}$. 
  Let $\Delta: (1)_{n \in \mathbb{N}} \to (\mathbb{N})_{n \in \mathbb{N}}$ be the map that at level $n$ 
  is just the map $\{n\}: 1 \to \mathbb{N}$. We can think of this map as a ``diagonal map". 
  Let $(P^m_n)_{n \in \mathbb{N}}$ be the following pullbacks 
 
 \begin{center}
  \pbsq{(P^m_n)_{n \in \mathbb{N}}}{(1)_{n \in \mathbb{N}}}{(1)_{n \in \mathbb{N}}}{(\mathbb{N})_{n \in \mathbb{N}}}{
  \varphi_m}{}{\Delta}{ \{ m \} }
 \end{center}
 
 By descent if the cocone formed by the maps $\{ \{ m \}: (1)_{n \in \mathbb{N}} \to (\mathbb{N})_{n \in \mathbb{N}} \}_{m \in \mathbb{N}}$ 
 is a colimiting 
 cocone, then the cocone formed by $\{\varphi_m: (P^m_n)_{n \in \mathbb{N}} \to (1)_{n \in \mathbb{N}} \}_{m \in \mathbb{N}}$ is also a colimiting cocone. 
 However, the maps $\Delta,\{m\} : (1)_{n \in \mathbb{N}} \to (\mathbb{N})_{n \in \mathbb{N}}$ only coincide when $n = m$ 
 and disagree otherwise. Thus for $n > m$, we have $P^m_n = \emptyset$, which implies that $(1)_{n \in \mathbb{N}}$ is isomorphic 
 to the colimit $\coprod_{n \in \mathbb{N}} \emptyset \cong \emptyset$, which only happens in the topos with one element.
 \par 
 Notice, by the explanation above, that $\prod_\Phi \Kan$ also does not have infinite colimits.
 Thus we cannot use traditional methods to study its homotopy theory (such as define truncations and prove Blakers-Massey theorem).
 This necessitates developing elementary $(\infty,1)$-topos theory and proving those results in the elementary context.
 For such an elementary approach and a careful analysis of truncations in $\prod_\Phi \Kan$ see \cite{rasekh2018truncations}.
 
 \begin{remone}
  For an alternative argument why $\prod_\mathbb{N} \Kan$ does not have infinite limits see \cite[Example 3.19]{bss2020chromatic}.
 \end{remone}

 \section{Future Directions} \label{Sec:Future}
  Up to here we answered the important question on whether there is an elementary $(\infty,1)$-topos that is not a 
  Grothendieck $(\infty,1)$-topos. 
  In this section we want to pose some new questions about elementary $(\infty,1)$-topos theory motivated by the response. 
 
 {\bf Homotopy Filter Quotient:} 
 Although the filter quotient construction enables us to construct new elementary $(\infty,1)$-toposes,
 it does have certain limitations. In particular, a filter on $\E$ is completely determined by a filter on the underlying elementary 
 topos $\tau_0\E$ as we are just using $\Sub(1)$, which is a subcategory of both categories.
 This is particularly problematic when we have interesting categorical data as we will illustrate in the following example.
 
 \begin{exone}
  Let $\Kan$ be the Kan enriched category of Kan complexes and let $K$ be a connected Kan complex. Then the slice category 
  $\Kan_{/K}$ is also an elementary $(\infty,1)$-topos with final object $id_K: K \to K$, which has only two subobjects, namely 
  $\emptyset \to K$ and $id_K : K \to K$. Thus $\Kan_{/K}$ only has trivial filter quotients (see \cref{Ex:Triv Filter}).
 \end{exone}

 Ideally we would like to have a notion of filter quotient which allows us to take the higher homotopical data of $K$ into account 
 and use that to form a quotient, which should lead us to a notion of {\it homotopy filter quotient}.
 
 {\bf Filter Quotients of Sheaves:}
 We give explicit examples of filter quotients for the simplest of all Grothendieck $(\infty,1)$-toposes, namely $\Kan^S$, 
 where $S$ is a set (\cref{Ex:KanS Phi not GHT}). 
 However, given that the filter construction works for all Grothendieck $(\infty,1)$-toposes the next step is to 
 construct the filter quotient on the category of sheaves. 
 
 {\bf Non-standard Models of Spaces:}
 One long standing goal of elementary $(\infty,1)$-topos theory is to study models for spaces (similar to how 
 elementary toposes are used to study models of set theory) and the filter quotient might be a step in developing such models. 
 Concretely, let $S$ be an infinite set and $\Phi$ a non-principal ultrafilter (\cref{Def:Ultrafilter}) on the power set $P(S)$. 
 Then we can construct the ultra product $\prod_\Phi \Kan$. The underlying elementary topos is the ultra product 
 $\prod_\Phi \set$, which has many interesting properties \cite[A2.2]{johnstone2002elephanti}:
 \begin{enumerate}
  \item It is Boolean.
  \item It is generated by the final object.
  \item It doesn't have infinite colimits.
 \end{enumerate}
 Thus $\prod_\Phi \set$ shares many properties with the category of sets. This suggests that $\prod_\Phi \Kan$ should behave similar to the 
 $(\infty,1)$-category of Kan complexes. An important first step would be to prove that $\prod_\Phi\Kan$ is generated by the final object.
 
 {\bf Models of Homotopy Type Theory:}
 One important goal of elementary $(\infty,1)$-topos theory is to construct models of homotopy type theory.
 We already know that every Grothendieck $(\infty,1)$-topos is a model of homotopy type theory \cite{shulman2019inftytoposunivalent}.
 However, as of now there are no other known models. 
 Given that we have an explicit construction of the filter quotient the hope is that we can show the filter quotient 
 of a Grothendieck $(\infty,1)$-topos is also a model of homotopy type theory.
 
 \bibliographystyle{alpha}
 \bibliography{main}

\begin{thebibliography}{EKMM95}

\bibitem[ABFJ18]{abfj2018goodwillie}
M.~Anel, G.~Biedermann, E.~Finster, and A.~Joyal.
\newblock Goodwillie's calculus of functors and higher topos theory.
\newblock {\em J. Topol.}, 11(4):1100--1132, 2018.

\bibitem[ABFJ20]{abfj2017blakersmassey}
Mathieu Anel, Georg Biedermann, Eric Finster, and Andr\'{e} Joyal.
\newblock A generalized {B}lakers-{M}assey theorem.
\newblock {\em J. Topol.}, 13(4):1521--1553, 2020.

\bibitem[Ada95]{adams1995stable}
J.~F. Adams.
\newblock {\em Stable homotopy and generalised homology}.
\newblock Chicago Lectures in Mathematics. University of Chicago Press,
  Chicago, IL, 1995.
\newblock Reprint of the 1974 original.

\bibitem[AJ82]{adelmanjohnstone1982serreclasses}
M.~Adelman and P.~T. Johnstone.
\newblock Serre classes for toposes.
\newblock {\em Bull. Austral. Math. Soc.}, 25(1):103--115, 1982.

\bibitem[AK65a]{axkochen1965diophantinei}
James Ax and Simon Kochen.
\newblock Diophantine problems over local fields. {I}.
\newblock {\em Amer. J. Math.}, 87:605--630, 1965.

\bibitem[AK65b]{axkochen1965diophantineii}
James Ax and Simon Kochen.
\newblock Diophantine problems over local fields. {II}. {A} complete set of
  axioms for {$p$}-adic number theory.
\newblock {\em Amer. J. Math.}, 87:631--648, 1965.

\bibitem[AK66]{axkochen1966diophantineiii}
James Ax and Simon Kochen.
\newblock Diophantine problems over local fields. {III}. {D}ecidable fields.
\newblock {\em Ann. of Math. (2)}, 83:437--456, 1966.

\bibitem[Ber07]{bergner2007bergnermodelcat}
Julia~E. Bergner.
\newblock A model category structure on the category of simplicial categories.
\newblock {\em Trans. Amer. Math. Soc.}, 359(5):2043--2058, 2007.

\bibitem[Ber10]{bergner2010survey}
Julia~E. Bergner.
\newblock A survey of {$(\infty,1)$}-categories.
\newblock In {\em Towards higher categories}, volume 152 of {\em IMA Vol. Math.
  Appl.}, pages 69--83. Springer, New York, 2010.

\bibitem[Ber18]{bergner2018book}
Julia~E. Bergner.
\newblock {\em The homotopy theory of {$(\infty, 1)$}-categories}, volume~90 of
  {\em London Mathematical Society Student Texts}.
\newblock Cambridge University Press, Cambridge, 2018.

\bibitem[Bro27]{brouwer1927allcontinuous}
L.~E.~J. Brouwer.
\newblock \"{U}ber {D}efinitionsbereiche von- {F}unktionen.
\newblock {\em Math. Ann.}, 97(1):60--75, 1927.

\bibitem[BSS20]{bss2020chromatic}
Tobias Barthel, Tomer~M. Schlank, and Nathaniel Stapleton.
\newblock Chromatic homotopy theory is asymptotically algebraic.
\newblock {\em Invent. Math.}, 220(3):737--845, 2020.

\bibitem[BV73]{boardmanvogt1973qcats}
J.~M. Boardman and R.~M. Vogt.
\newblock {\em Homotopy invariant algebraic structures on topological spaces}.
\newblock Lecture Notes in Mathematics, Vol. 347. Springer-Verlag, Berlin-New
  York, 1973.

\bibitem[Car37a]{cartan1937ultrafiltres}
Henri Cartan.
\newblock Filtres et ultrafiltres.
\newblock {\em Comptes Rendus Hebdomadaires des Seances de l'Academie des
  Sciences, Paris}, 205:777--779, 1937.

\bibitem[Car37b]{cartan1937filters}
Henri Cartan.
\newblock Th{\'e}orie des filtres.
\newblock {\em CR Acad. Sci. Paris}, 205:595--598, 1937.

\bibitem[Chu40]{church1940typetheory}
Alonzo Church.
\newblock A formulation of the simple theory of types.
\newblock {\em J. Symbolic Logic}, 5:56--68, 1940.

\bibitem[CK90]{changkeisler1990modeltheory}
C.~C. Chang and H.~J. Keisler.
\newblock {\em Model theory}, volume~73 of {\em Studies in Logic and the
  Foundations of Mathematics}.
\newblock North-Holland Publishing Co., Amsterdam, third edition, 1990.

\bibitem[Coh63]{cohen1963forcing}
Paul Cohen.
\newblock The independence of the continuum hypothesis.
\newblock {\em Proc. Nat. Acad. Sci. U.S.A.}, 50:1143--1148, 1963.

\bibitem[DK80a]{dwyerkan1980calculatingsimplocalizations}
W.~G. Dwyer and D.~M. Kan.
\newblock Calculating simplicial localizations.
\newblock {\em J. Pure Appl. Algebra}, 18(1):17--35, 1980.

\bibitem[DK80b]{dwyerkan1980simplocalization}
W.~G. Dwyer and D.~M. Kan.
\newblock Simplicial localizations of categories.
\newblock {\em J. Pure Appl. Algebra}, 17(3):267--284, 1980.

\bibitem[EKMM95]{ekmm1995stable}
A.~D. Elmendorf, I.~K\v{r}\'{\i}\v{z}, M.~A. Mandell, and J.~P. May.
\newblock Modern foundations for stable homotopy theory.
\newblock In {\em Handbook of algebraic topology}, pages 213--253.
  North-Holland, Amsterdam, 1995.

\bibitem[EM45]{eilenbergmaclane1945categories}
Samuel Eilenberg and Saunders MacLane.
\newblock General theory of natural equivalences.
\newblock {\em Trans. Amer. Math. Soc.}, 58:231--294, 1945.

\bibitem[ES45]{eilenbergsteenrod1945axioms}
Samuel Eilenberg and Norman~E. Steenrod.
\newblock Axiomatic approach to homology theory.
\newblock {\em Proc. Nat. Acad. Sci. U.S.A.}, 31:117--120, 1945.

\bibitem[EZ50]{eilenbergzilber1950simplicial}
Samuel Eilenberg and J.~A. Zilber.
\newblock Semi-simplicial complexes and singular homology.
\newblock {\em Ann. of Math. (2)}, 51:499--513, 1950.

\bibitem[FMS63]{fraynemorelscott1962ultraproduct}
T.~Frayne, A.~C. Morel, and D.~S. Scott.
\newblock Reduced direct products.
\newblock {\em Fund. Math.}, 51:195--228, 1962/63.

\bibitem[GJ99]{goerssjardine1999simplicialhomotopytheory}
Paul~G. Goerss and John~F. Jardine.
\newblock {\em Simplicial homotopy theory}, volume 174 of {\em Progress in
  Mathematics}.
\newblock Birkh\"auser Verlag, Basel, 1999.

\bibitem[GK17]{gepnerkock2017univalence}
David Gepner and Joachim Kock.
\newblock Univalence in locally cartesian closed {$\infty$}-categories.
\newblock {\em Forum Math.}, 29(3):617--652, 2017.

\bibitem[G{\"{o}}d30]{godel1930completeness}
Kurt G{\"{o}}del.
\newblock Die {V}ollst\"{a}ndigkeit der {A}xiome des logischen
  {F}unktionenkalk\"{u}ls.
\newblock {\em Monatsh. Math. Phys.}, 37(1):349--360, 1930.

\bibitem[GSS19]{gss2019constructive}
Nicola Gambino, Christian Sattler, and Karol Szumi{\l}o.
\newblock The constructive {K}an-{Q}uillen model structure: two new proofs.
\newblock {\em arXiv preprint}, 2019.
\newblock \href{https://arxiv.org/abs/1907.05394v1}{arXiv:1907.05394v1}.

\bibitem[Hen19]{henry2019constructive}
Simon Henry.
\newblock A constructive account of the {K}an-{Q}uillen model structure and of
  {K}an's {E}x$^\infty$ functor.
\newblock {\em arXiv preprint}, 2019.
\newblock \href{https://arxiv.org/abs/1905.06160v1}{arXiv:1905.06160v1}.

\bibitem[Hir03]{hirschhorn2003modelcategories}
Philip~S. Hirschhorn.
\newblock {\em Model categories and their localizations}, volume~99 of {\em
  Mathematical Surveys and Monographs}.
\newblock American Mathematical Society, Providence, RI, 2003.

\bibitem[Hop31]{hopf1931fibration}
Heinz Hopf.
\newblock \"{U}ber die {A}bbildungen der dreidimensionalen {S}ph\"{a}re auf die
  {K}ugelfl\"{a}che.
\newblock {\em Math. Ann.}, 104(1):637--665, 1931.

\bibitem[Hov99]{hovey1999modelcategories}
Mark Hovey.
\newblock {\em Model categories}, volume~63 of {\em Mathematical Surveys and
  Monographs}.
\newblock American Mathematical Society, Providence, RI, 1999.

\bibitem[HSS00]{hss2000symmetricspectra}
Mark Hovey, Brooke Shipley, and Jeff Smith.
\newblock Symmetric spectra.
\newblock {\em J. Amer. Math. Soc.}, 13(1):149--208, 2000.

\bibitem[Isa19]{isaksen2019stablecomputation}
Daniel~C. Isaksen.
\newblock Stable stems.
\newblock {\em Mem. Amer. Math. Soc.}, 262(1269):viii+159, 2019.

\bibitem[Joh02a]{johnstone2002elephanti}
Peter~T. Johnstone.
\newblock {\em Sketches of an elephant: a topos theory compendium. {V}ol. 1},
  volume~43 of {\em Oxford Logic Guides}.
\newblock The Clarendon Press, Oxford University Press, New York, 2002.

\bibitem[Joh02b]{johnstone2002elephantsii}
Peter~T. Johnstone.
\newblock {\em Sketches of an elephant: a topos theory compendium. {V}ol. 2},
  volume~44 of {\em Oxford Logic Guides}.
\newblock The Clarendon Press, Oxford University Press, Oxford, 2002.

\bibitem[Joy08]{joyal2008notes}
Andr{\'e} Joyal.
\newblock Notes on quasi-categories, 2008.
\newblock \href{https://www.math.uchicago.edu/~may/IMA/Joyal.pdf}{Unpublished
  notes} (accessed 01.04.2021).

\bibitem[JT07]{joyaltierney2007qcatvssegal}
Andr\'{e} Joyal and Myles Tierney.
\newblock Quasi-categories vs {S}egal spaces.
\newblock In {\em Categories in algebra, geometry and mathematical physics},
  volume 431 of {\em Contemp. Math.}, pages 277--326. Amer. Math. Soc.,
  Providence, RI, 2007.

\bibitem[Kan57]{kan1957complexes}
Daniel~M. Kan.
\newblock On c. s. s. complexes.
\newblock {\em Amer. J. Math.}, 79:449--476, 1957.

\bibitem[Kan58]{kan1958homotopygroups}
Daniel~M. Kan.
\newblock A combinatorial definition of homotopy groups.
\newblock {\em Ann. of Math. (2)}, 67:282--312, 1958.

\bibitem[KL12]{kapulkinlumsdaine2012kanunivalent}
Krzysztof Kapulkin and Peter~LeFanu Lumsdaine.
\newblock The simplicial model of univalent foundations (after {V}oevodsky).
\newblock {\em arXiv preprint}, 2012.
\newblock \href{https://arxiv.org/abs/1211.2851v5}{arXiv:1211.2851v5}.

\bibitem[Law64]{lawvere1964elementarysets}
F.~William Lawvere.
\newblock An elementary theory of the category of sets.
\newblock {\em Proc. Nat. Acad. Sci. U.S.A.}, 52:1506--1511, 1964.

\bibitem[Lim59]{lima1959spectrum}
Elon~L. Lima.
\newblock The {S}panier-{W}hitehead duality in new homotopy categories.
\newblock {\em Summa Brasil. Math.}, 4:91--148 (1959), 1959.

\bibitem[LM21]{monaco2021eht}
Giulio Lo~Monaco.
\newblock Dependent products and 1-inaccessible universes.
\newblock {\em Theory Appl. Categ.}, 37:Paper No. 5, pp 107--143., 2021.

\bibitem[{\L}o{\'s}55]{los1955ultraproduct}
Jerzy {\L}o{\'s}.
\newblock Quelques remarques, th\'{e}or\`emes et probl\`emes sur les classes
  d\'{e}finissables d'alg\`ebres.
\newblock In {\em Mathematical interpretation of formal systems}, pages
  98--113. North-Holland Publishing Co., Amsterdam, 1955.

\bibitem[LS88]{lambekscott1988higherorderlogic}
J.~Lambek and P.~J. Scott.
\newblock {\em Introduction to higher order categorical logic}, volume~7 of
  {\em Cambridge Studies in Advanced Mathematics}.
\newblock Cambridge University Press, Cambridge, 1988.
\newblock Reprint of the 1986 original.

\bibitem[Lur09]{lurie2009htt}
Jacob Lurie.
\newblock {\em Higher topos theory}, volume 170 of {\em Annals of Mathematics
  Studies}.
\newblock Princeton University Press, Princeton, NJ, 2009.

\bibitem[Lur17]{lurie2017ha}
Jacob Lurie.
\newblock Higher algebra.
\newblock \href{http://www.math.ias.edu/~lurie/papers/HA.pdf}{Unpublished book}
  (accessed 01.04.2021), September 2017.

\bibitem[Lur18]{lurie2018sag}
Jacob Lurie.
\newblock Spectral algebraic geometry.
\newblock
  \href{https://www.math.ias.edu/~lurie/papers/SAG-rootfile.pdf}{Unpublished
  notes} (accessed 01.04.2021), February 2018.

\bibitem[ML75a]{martinlof1975inttypetheories}
Per Martin-L\"{o}f.
\newblock About models for intuitionistic type theories and the notion of
  definitional equality.
\newblock In {\em Proceedings of the {T}hird {S}candinavian {L}ogic {S}ymposium
  ({U}niv. {U}ppsala, {U}ppsala, 1973)}, pages 81--109. Stud. Logic Found.
  Math., Vol. 82, 1975.

\bibitem[ML75b]{martinlof1975inttypetheoriespredicative}
Per Martin-L\"{o}f.
\newblock An intuitionistic theory of types: predicative part.
\newblock In {\em Logic {C}olloquium '73 ({B}ristol, 1973)}, pages 73--118.
  Studies in Logic and the Foundations of Mathematics, Vol. 80. 1975.

\bibitem[ML98]{maclane1998categories}
Saunders Mac~Lane.
\newblock {\em Categories for the working mathematician}, volume~5 of {\em
  Graduate Texts in Mathematics}.
\newblock Springer-Verlag, New York, second edition, 1998.

\bibitem[MLM94]{maclanemoerdijk1994topos}
Saunders Mac~Lane and Ieke Moerdijk.
\newblock {\em Sheaves in geometry and logic}.
\newblock Universitext. Springer-Verlag, New York, 1994.
\newblock A first introduction to topos theory, Corrected reprint of the 1992
  edition.

\bibitem[MM02]{mandellmay2002orthogonal}
M.~A. Mandell and J.~P. May.
\newblock Equivariant orthogonal spectra and {$S$}-modules.
\newblock {\em Mem. Amer. Math. Soc.}, 159(755):x+108, 2002.

\bibitem[Pal97]{palmgren1997nonstandardanalysis}
Erik Palmgren.
\newblock A sheaf-theoretic foundation for nonstandard analysis.
\newblock {\em Ann. Pure Appl. Logic}, 85(1):69--86, 1997.

\bibitem[Poi95]{poincare1895analysissitus}
Henri Poincar{\'e}.
\newblock Analysis situs.
\newblock 1895.

\bibitem[Qui67]{quillen1967modelcats}
Daniel~G. Quillen.
\newblock {\em Homotopical algebra}.
\newblock Lecture Notes in Mathematics, No. 43. Springer-Verlag, Berlin-New
  York, 1967.

\bibitem[Ras17]{rasekh2017left}
Nima Rasekh.
\newblock Yoneda lemma for simplicial spaces.
\newblock {\em arXiv preprint}, 2017.
\newblock \href{https://arxiv.org/abs/1711.03160v3}{arXiv:1711.03160v3}.

\bibitem[Ras18a]{rasekh2018truncations}
Nima Rasekh.
\newblock An elementary approach to truncations.
\newblock {\em arXiv preprint}, 2018.
\newblock \href{https://arxiv.org/abs/1812.10527v2}{arXiv:1812.10527v2}.

\bibitem[Ras18b]{rasekh2018elementarytopos}
Nima Rasekh.
\newblock A theory of elementary higher toposes.
\newblock {\em arXiv preprint}, 2018.
\newblock \href{https://arxiv.org/abs/1805.03805v2}{arXiv:1805.03805v2}.

\bibitem[Ras21]{rasekh2021nno}
Nima Rasekh.
\newblock Every elementary higher topos has a natural number object.
\newblock {\em Theory Appl. Categ.}, 37:Paper No. 13, pp 337--377, 2021.

\bibitem[Rez01]{rezk2001css}
Charles Rezk.
\newblock A model for the homotopy theory of homotopy theory.
\newblock {\em Trans. Amer. Math. Soc.}, 353(3):973--1007, 2001.

\bibitem[Rez10a]{rezk2010thetanspaces}
Charles Rezk.
\newblock A {C}artesian presentation of weak {$n$}-categories.
\newblock {\em Geom. Topol.}, 14(1):521--571, 2010.

\bibitem[Rez10b]{rezk2010toposes}
Charles Rezk.
\newblock Toposes and homotopy toposes (version 0.15).
\newblock
  \href{https://faculty.math.illinois.edu/~rezk/homotopy-topos-sketch.pdf}{Unpublished
  notes} (accessed 01.04.2021), 2010.

\bibitem[RS17]{riehlshulman2017rezktypes}
Emily Riehl and Michael Shulman.
\newblock A type theory for synthetic {$\infty$}-categories.
\newblock {\em High. Struct.}, 1(1):147--224, 2017.

\bibitem[Rus08]{russell1908typetheory}
Bertrand Russell.
\newblock Mathematical {L}ogic as {B}ased on the {T}heory of {T}ypes.
\newblock {\em Amer. J. Math.}, 30(3):222--262, 1908.

\bibitem[RV17]{riehlverity2017inftycosmos}
Emily Riehl and Dominic Verity.
\newblock Fibrations and {Y}oneda's lemma in an {$\infty$}-cosmos.
\newblock {\em J. Pure Appl. Algebra}, 221(3):499--564, 2017.

\bibitem[sga72a]{sga1972tome1}
{\em Th\'{e}orie des topos et cohomologie \'{e}tale des sch\'{e}mas. {T}ome 1:
  {T}h\'{e}orie des topos}.
\newblock Lecture Notes in Mathematics, Vol. 269. Springer-Verlag, Berlin-New
  York, 1972.
\newblock S\'{e}minaire de G\'{e}om\'{e}trie Alg\'{e}brique du Bois-Marie
  1963--1964 (SGA 4), Dirig\'{e} par M. Artin, A. Grothendieck, et J. L.
  Verdier. Avec la collaboration de N. Bourbaki, P. Deligne et B. Saint-Donat.

\bibitem[sga72b]{sga1972tome2}
{\em Th\'{e}orie des topos et cohomologie \'{e}tale des sch\'{e}mas. {T}ome 2}.
\newblock Lecture Notes in Mathematics, Vol. 270. Springer-Verlag, Berlin-New
  York, 1972.
\newblock S\'{e}minaire de G\'{e}om\'{e}trie Alg\'{e}brique du Bois-Marie
  1963--1964 (SGA 4), Dirig\'{e} par M. Artin, A. Grothendieck et J. L.
  Verdier. Avec la collaboration de N. Bourbaki, P. Deligne et B. Saint-Donat.

\bibitem[sga73]{sga1972tome3}
{\em Th\'{e}orie des topos et cohomologie \'{e}tale des sch\'{e}mas. {T}ome 3}.
\newblock Lecture Notes in Mathematics, Vol. 305. Springer-Verlag, Berlin-New
  York, 1973.
\newblock S\'{e}minaire de G\'{e}om\'{e}trie Alg\'{e}brique du Bois-Marie
  1963--1964 (SGA 4), Dirig\'{e} par M. Artin, A. Grothendieck et J. L.
  Verdier. Avec la collaboration de P. Deligne et B. Saint-Donat.

\bibitem[Shu19]{shulman2019inftytoposunivalent}
Michael Shulman.
\newblock All ($\infty$, 1)-toposes have strict univalent universes.
\newblock {\em arXiv preprint}, 2019.
\newblock \href{https://arxiv.org/abs/1904.07004v2}{arXiv:1904.07004v2}.

\bibitem[Sko34]{skolem1934earlyultraproducts}
Th. Skolem.
\newblock {\"U}ber die nicht-charakterisierbarkeit der {Z}ahlenreihe mittels
  endlich oder abz{\"a}hlbar unendlich vieler {A}ussagen mit ausschliesslich
  {Z}ahlenvariablen.
\newblock {\em Fundamenta mathematicae}, 23(1):150--161, 1934.

\bibitem[Sko55]{skolem1955models}
Th. Skolem.
\newblock Peano's axioms and models of arithmetic.
\newblock In {\em Mathematical interpretation of formal systems}, pages 1--14.
  North-Holland Publishing Co., Amsterdam, 1955.

\bibitem[SV00]{suslinvoevodsky2000blochkato}
Andrei Suslin and Vladimir Voevodsky.
\newblock Bloch-{K}ato conjecture and motivic cohomology with finite
  coefficients.
\newblock In {\em The arithmetic and geometry of algebraic cycles ({B}anff,
  {AB}, 1998)}, volume 548 of {\em NATO Sci. Ser. C Math. Phys. Sci.}, pages
  117--189. Kluwer Acad. Publ., Dordrecht, 2000.

\bibitem[Tar56]{tarski1956semantics}
Alfred Tarski.
\newblock {\em Logic, semantics, metamathematics. {P}apers from 1923 to 1938}.
\newblock Oxford at the Clarendon Press, 1956.
\newblock Translated by J. H. Woodger.

\bibitem[Tie72]{tierney1972elementarycontinuum}
Myles Tierney.
\newblock Sheaf theory and the continuum hypothesis.
\newblock In {\em Toposes, algebraic geometry and logic ({C}onf., {D}alhousie
  {U}niv., {H}alifax, {N}.{S}., 1971)}, pages 13--42. Lecture Notes in Math.,
  Vol. 274. 1972.

\bibitem[Tie73]{tierney1973elementarytopos}
M.~Tierney.
\newblock Axiomatic sheaf theory: some constructions and applications.
\newblock In {\em Categories and commutative algebra ({C}.{I}.{M}.{E}., {III}
  {C}iclo, {V}arenna, 1971)}, pages 249--326. 1973.

\bibitem[{Uni}13]{hottbook2013}
The {Univalent Foundations Program}.
\newblock {\em Homotopy Type Theory: Univalent Foundations of Mathematics}.
\newblock \url{https://homotopytypetheory.org/book}, Institute for Advanced
  Study, 2013.

\bibitem[Ver19]{vergura2019localization}
Marco Vergura.
\newblock Localization theory in an $\infty$-topos.
\newblock {\em arXiv preprint}, 2019.
\newblock \href{https://arxiv.org/abs/1907.03836v1}{arXiv:1907.03836v1}.

\bibitem[Voe98]{voevodsky1998motivicorigin}
Vladimir Voevodsky.
\newblock {$\bold A^1$}-homotopy theory.
\newblock In {\em Proceedings of the {I}nternational {C}ongress of
  {M}athematicians, {V}ol. {I} ({B}erlin, 1998)}, number Extra Vol. I, pages
  579--604, 1998.

\bibitem[Voe14]{voevodsky2014origins}
Vladimir Voevodsky.
\newblock The origins and motivations of univalent foundations.
\newblock {\em The Institute Letter}, pages 8--9, 2014.

\bibitem[Vol75]{volger1975filter}
Hugo Volger.
\newblock Ultrafilters, ultrapowers and finiteness in a topos.
\newblock {\em J. Pure Appl. Algebra}, 6(3):345--356, 1975.

\bibitem[Whi49a]{whitehead1949combinatoriali}
J.~H.~C. Whitehead.
\newblock Combinatorial homotopy. {I}.
\newblock {\em Bull. Amer. Math. Soc.}, 55:213--245, 1949.

\bibitem[Whi49b]{whitehead1949combinatorialii}
J.~H.~C. Whitehead.
\newblock Combinatorial homotopy. {II}.
\newblock {\em Bull. Amer. Math. Soc.}, 55:453--496, 1949.

\end{thebibliography}

\end{document}